\documentclass{amsart}
\usepackage{amsfonts}
\usepackage{amscd}
\usepackage{amsmath}

\theoremstyle{plain} \numberwithin{equation}{section}
\newtheorem{theorem}{Theorem}[section]
\newtheorem{corollary}[theorem]{Corollary}
\newtheorem{conjecture}{Conjecture}
\newtheorem{lemma}[theorem]{Lemma}
\newtheorem{proposition}[theorem]{Proposition}
\theoremstyle{definition}
\newtheorem{definition}[theorem]{Definition}

\newtheorem{remark}[theorem]{Remark}
\newtheorem{example}[theorem]{Example}
\newtheorem{question}{Question}
\numberwithin{equation}{section}

\newcommand{\cD}{\mathcal{D}}

\newcommand{\cB}{\mathcal{B}}

\newcommand{\cO}{\mathcal{O}}
\newcommand{\cH}{\mathcal{H}}

\newcommand{\cK}{\mathcal{K}}
\newcommand{\cT}{\mathcal{T}}

\newcommand{\GL}{\mathrm{GL}}

\newcommand{\SO}{\mathrm{SO}}

\newcommand{\fa}{\mathfrak{a}}

\newcommand{\fc}{\mathfrak{c}}

\newcommand{\gO}{\Omega}
\newcommand{\cL}{\mathcal {L}}

\newcommand{\F}{\mathbb{F}}
\newcommand{\R}{\mathbb{R}}
\newcommand{\C}{\mathbb{C}}
\newcommand{\J}{\mathbb{J}}
\newcommand{\N}{\mathbb{N}}
\newcommand{\Z}{\mathbb{Z}}

\newcommand{\brn}{\R^n}
\newcommand{\pr}{\mathop{\mathrm{pr}} }

\newcommand{\id}{\mathrm{id}}

\newcommand{\abs}[1]{\left\vert#1\right\vert}

\theoremstyle{plain}

\begin{document}
\title[Groups, Wavelets, and Wavelet Sets]{Groups, Wavelets, and Wavelet Sets}
\author{Gestur \'{O}lafsson and Darrin Speegle}
\address{Department of Mathematics, Louisiana State
University, Baton Rouge,
LA 70803, USA} \email{olafsson@math.lsu.edu}
\address{Department of Mathematics, Saint Louis University,
St. Louis, MO
63103} \email{speegled@slu.edu} 
\thanks{Research supported by NSF grants DMS-0070607 and
DMS-0139783}

\begin{abstract}
Wavelet and frames have become a widely used tool in
mathematics,
physics, and applied science during the last decade. This
article
gives an overview over some well known results about the
continuous and discrete wavelet transforms and groups acting
on
$\mathbb{R}^n$.  We also show how this action can give rise
to
wavelets, and in particular, MSF wavelets)in
$L^2(\mathbb{R}^n)$.
\end{abstract}

\maketitle

\section*{Introduction}

\noindent
The classical wavelet system consists of a single function
$\psi\in L^2(\R)$ such that $\{2^{j/2}\psi(2^jx + k) \mid j,k
\in
\Z\}$ is an orthonormal basis for $L^2(\R)$.  There has
been quite a bit of recent interest in relaxing various aspects
of the definition of wavelets, in particular
in higher dimensions.  For example, one can allow multiple
functions ${\psi^i, \ldots, \psi^L}$, an arbitrary matrix of
dilations, and an arbitrary lattice of translations.  One
could relax even further to allow a group of dilations,
or perhaps even
just a set of dilations and translations.  A first question
one would ask, then, is: for which collections of dilations and
translations do there exist wavelets?  We will begin by
reviewing some well-known results concerning this central question.
Then, we will show that there is a fundamental connection between
the papers of Dai, Diao, Gu and Han \cite{DDGH}, Fabec and
{\'O}lafsson
\cite{FO02}, Laugesen, Weaver, Weiss and Wilson
\cite{LWWW2002}, and Wang \cite {W02}. One argument that a survey
paper such as this one is usefull is that, even though these eleven authors
are active in the field, there is only one cross-reference of
the above papers in the references of the other papers.

We now describe briefly the connection betwen the papers
listed above.  All four papers are concerned with constructing
reproducing systems consisting of dilations and translations
of a function. That is, they consider triples $(\cD,\cT, M)$,
where $\cD$ is some collection of invertible matrices, $\cT$ is
some collection of points in $\brn$, and $M$ a non-trivial closed subspace of
$L^2(\brn)$. Then, they ask whether there is a function $\psi$ such that
$\{\psi_{a,k}=\abs{\det a}^{1/2}\psi(a(x)+k)\mid a\in \cD, k\in
\cT\}$ is a frame, normalized tight frame, or even a orthonormal basis for $M$

In \cite{DDGH}, it is assumed that $\cD = \{a^j \mid j\in \Z\}$
for some expansive matrix $a $, that $\cT = \Z^n$, and that $M$ is an
$a$-invariant subspace of $L^2(\brn)$.  In \cite{FO02}, the
assumptions are that,  $\cD$ is constructed as a subset of a
particular type of group $H$, that $\cT$ is a full
rank lattice depending on $H$,
and finally that $M$ is of the form
$M=\{f\in L^2(\R^n)\mid \textrm{Supp}(\hat{f})\subseteq \overline{\cO}\}$,
where $\cO\subset \R^n$ is an open $H$-orbit.
In \cite{LWWW2002}, it is
assumed that $\cD$ is a group, $\cT=\Z^n$, and
$M=L^2(\brn)$.  In \cite{W02}, it is assumed that $\cD$ and $\cT$ satisfy
non-algebraic conditions relating to the existence of
fundamental regions (see Section \ref{s-wset} for details) and $M= L^2(\brn)$.
Moreover, all four papers - either explicitly or
implicitly - are concerned primarily with the existence of
functions of the form $\hat\psi= \chi_\Omega$.

When put in this general framework, it becomes clear that
the four papers are related in spirit and scope.  What we will show
below is that they are also related in that results in \cite{DDGH}
can be used to remove technical assumptions from results in
\cite{W02}.  The improved results in \cite{W02} can then be
used to improve the results in \cite{LWWW2002} and \cite{FO02}.
We will improve the results in \cite{FO02} by removing the
dependence of the lattice on the group, and by constructing an
orthonormal basis where a normalized tight frame was constructed before.
The proof of the main Theorem in \cite{FO02} will also be
simplified.
Finally, we improve the results in \cite{LWWW2002} by
replacing normalized tight frame system with a wavelet system.

We will attempt to make these technical improvements to the
theorems in these papers with a minimal amount of technical
work. In particular, where possible, we will apply theorem quoting
proofs.  The primary exception to this is Theorem
\ref{dlsimprovement}, where we essentially need to check
that the details of an argument in \cite{DLSII} go through in a
slightly more general setting.

\section{Wavelet sets}\label{s-wset}

\noindent We start this section by recalling some simple
definitions and facts about wavelets, wavelet sets, and
tilings.
For a measureable set $\Omega\subseteq \R^n$ we denote by
$\chi_\Omega$ the indicator function of the set $\Omega$ and
by
$|\Omega |=\int \chi_\Omega (x)\, dx$ the measure of
$\Omega$.

\begin{definition}
A countable collection $\{\Omega_j\}$ of subsets of
$\mathbb{R}^n$
is a (measurable) tiling of $\mathbb{R}^n$ if
$|\mathbb{R}^n\setminus \bigcup_j\Omega_j|=0$, and
$|\Omega_i\cap
\Omega_j|=0$ for $i\not=j$.
\end{definition}

\begin{definition}
Let $\mathcal{T}\subset \mathbb{R}^n$ and
$\mathcal{D}\subset \mathrm{GL} (n,%
\mathbb{R})$. We say that $\mathcal{D}$ is a
\textit{multiplicative tiling
set} of $\mathbb{R}^n$ if there exists a set $\Omega\subset
\mathbb{R%
}^n$ of positive measure such that $\{d(\Omega)\mid d\in
\mathcal{D}\}$ is a tiling of $\mathbb{R}^n$. The set
$\Omega$ is
said to be a \textit{multiplicative}
$\mathcal{D}$-\textit{tile}.
We say $\mathcal{D}$ is a \textit{bounded multiplicative
tiling
set} of $\mathbb{R}^n$ if there is a multiplicative
$\mathcal{D}$-tile $\Omega$ which is bounded and such that
$0\not
\in \overline{\Omega}$.

Similarly, we say that $\mathcal{T}$ is an \textit{additive
tiling
set} of $\brn$ if there exists a set $\Omega\subset \brn$
such
that $\{\Omega + x\mid x\in \cT\}$ is a tiling of $\brn$. 
The set
$\Omega$ is said to be an \textit{(additive)} $\cT$
\textit{tile}.
Again, we add the word bounded if $\Omega$ can be chosen to
be a
bounded set (with no restriction on being bounded away from
0).

A set $\Omega$ is a \textit{$(\cD,\cT)$ tiling set} if it is
a
$\cD$ multiplicative tiling set and a $\cT$ additive tiling
set.

\end{definition}

Note that this definition does not coincide with the
definition of
Wang \cite{W02}.  Wang defines a multiplicative tiling set
to be
what we have defined to be a bounded multiplicative tiling
set.
We feel that boundedness properties of $\mathcal{D}$-tiles
are
interesting properties, but they should not be part of a
definition of tiling.

Multiplicative and additive tilings of $\mathbb{R}^n$ show
up in
wavelet theory and other branches of analysis in a natural
way.

\begin{definition}
A function $\varphi\in L^2(\mathbb{R}^n)$ is called a
\textit{wavelet} if there exists a subset
$\mathcal{D}\subset
\mathrm{GL} (n,\mathbb{R})$ and a subset $\mathcal{T}\subset
\mathbb{R}^n$ such that
\begin{equation*}
\mathcal{W}(\varphi;\mathcal{D},\mathcal{T}):=\{|\det
d|^{1/2}\varphi(dx+k)\mid d\in \mathcal{D},\, k
\in\mathcal{T}\}
\end{equation*}
forms an orthonormal basis for $L^2(\mathbb{R}^n)$. The set
$\mathcal{D}$ is called the \textit{dilation set} for
$\varphi$,
the set $\mathcal{T}$ is called the \textit{translation set}
for
$\varphi$, and we say that $\varphi$ is a
$\mathbf{(\mathcal{D},\mathcal{T})}$\textit{-wavelet}.
\end{definition}

Normalize the Fourier transform by
\begin{equation*}
\mathcal{F}(f)(\lambda )=\hat{f}(\lambda
)=\int_{\mathbb{R}^{n}}f(x)e^{2\pi i(\lambda ,x)}\,dx\,.
\end{equation*}%
We set $f^{\vee }(x)=\hat{f}(-x)$. Then $f=(\hat{f})^{\vee
}$. For
simplicity we set $e_{\lambda }(x)=e^{2\pi i(\lambda ,x)}$.

\begin{definition}
Let $\Omega \subset \mathbb{R}^{n}$ be measurable with
positive,
but finite measure. We say that $\Omega $ is \textit{a
wavelet
set} if there exists a
pair $(\mathcal{D},\mathcal{T})$, with $\mathcal{D}\subset
\mathrm{GL}(n,%
\mathbb{R})$ and $\mathcal{T}\subset \mathbb{R}^{n}$ such
that
$\chi
_{\Omega }^{\vee }$ is a
$(\mathcal{D},\mathcal{T})$-wavelet. If $\chi _{\gO%
}^{\vee }$ is a $(\mathcal{D},\mathcal{T})$-wavelet, then we
say that $%
\Omega $ is a
$\mathbf{(\mathcal{D},\mathcal{T})}$\textit{-wavelet
set}.
\end{definition}

\begin{definition}
A measurable set $\Omega \subset \mathbb{R}^{n}$ with finite
positive
measure is called \textit{a spectral set} if there exists a
set $\mathcal{T}%
\subset \mathbb{R}^{n}$ such that the sequence of functions
$\left\{ e_{\lambda }\right\} _{\lambda \in \mathcal{T}}$
forms an
orthogonal basis for $L^{2}(\Omega )$. If this is the case
we say
that $\mathcal{T}$ is the \textit{spectrum of }$\Omega $,
and say
that  $(\Omega ,\mathcal{T})$ is a spectral pair
\end{definition}

Now, after given this list of definitions, let us recall
some
results, questions, and conjectures on how these concepts
are tied
together.  A first result, which has appeared in several
places
\cite{DL,HKLS,HW} is

\begin{theorem}
A measurable set $\Omega\subset \mathbb{R}$ is a wavelet set
for the pair $%
\mathcal{D}=\{2^n\mid n\in \mathbb{Z}\}$ and
$\mathcal{T}=\mathbb{Z}$ if and only if $\Omega$ is a
$(\mathcal{D},\mathcal{T})$-tiling set.
\end{theorem}

For the general case we have now the following two related
questions:

\begin{question}[Wang, \protect\cite{W02}]
For which sets $\mathcal{D}\subset \GL(n, \mathbb{R})$,
$\mathcal{T}\subset \mathbb{R}^n$ do there exist
$(\mathcal{D},
\mathcal{T})$-wavelets?
\end{question}

\begin{question}
For which sets $\mathcal{D}\subset \GL(n, \mathbb{R})$,
$\mathcal{T}\subset \mathbb{R}^n$ do there exist
$(\mathcal{D},\mathcal{T})$-wavelet sets?
\end{question}

Clearly, if there exists a $(\mathcal{D}, \mathcal{T})$
wavelet
set, then there exists a $(\mathcal{D}, \mathcal{T})$
wavelet,
but, what is interesting, is that the converse may also be
true.
In particular, there are currently no examples known of sets
$(\mathcal{D}, \mathcal{T})$ for which there exist wavelets,
but
for which there do not exist $(\mathcal{D}, \mathcal{T})$
wavelet
sets. Therefore we can state the third natural question:

\begin{question}
\label{msfsame} Is it true that if there exists a
$(\mathcal{D},
\mathcal{T}) $ wavelet, then there exists a $(\mathcal{D},
\mathcal{T})$ wavelet set?
\end{question}

So far, all evidence points to a positive answer for
question
\ref{msfsame}. (Though, we should point out that question
\ref{msfsame} has mostly been thought about in the case that
$\mathcal{D}$ is a singly generated group and $\mathcal{T}$
is a
full rank lattice, so it is possible that there is a
relatively
easy counterexample to the question posed in this
generality.)
When $\mathcal{D}$ is generated by a single matrix $a$ and
$\mathcal{T}$ is a lattice, it is known \cite{DLS} that if
$a$ is
expansive, then there exist $(\mathcal{D}, \mathcal{T})$
wavelet
sets. Moreover, it is also known \cite {Bow, BS, CS} in the
expansive case that there exist $(\mathcal{D}, \mathcal{T})$
wavelets that do not come from a wavelet set if and only if
there is a $%
j\not= 0$ such that $(a^T)^j (\mathcal{T}^*) \cap \mathcal{T}^*
\not = \{0\}$. In particular, for \textit{most} pairs of
this
type, the \textit{only} wavelets that exist come from
wavelet
sets. When $\mathcal{D}$ is generated by a not necessarily
expansive matrix $a$ and $\mathcal{T}$ is a lattice, then
the
handful of $(\mathcal{D}, \mathcal{T})$ wavelets known all
come from $(\mathcal{D}, \mathcal{T})$ wavelet sets.

There is also a stronger version of question \ref{msfsame}
due to
Larson \cite{L} in the one dimensional case.

\begin{question}[Larson, \protect\cite{L}]
\label{larsonq} Is it true that if $\psi $ is a
$(\mathcal{D},\mathcal{T})$ wavelet, then there is a
$(\mathcal{D},\mathcal{T})$ wavelet set $E\subseteq
\mathrm{supp}(\hat{\psi})$?
\end{question}

This problem is open even for the \textquotedblleft
classical"
case of dimension 1 with dilations by powers of 2 and
translations
by integers. We name two partial answers. The first is given
by
Rzeszotnik in his PhD Thesis, and the second is due to
Rzeszotnik
and the second author of this paper.

\begin{theorem}[Rzeszotnik, \protect\cite{R} Corollary 3.10]
Every multiresulution analysis (MRA) $(2^j,\Z)$ wavelet
contains
in its supportthe support of its Fourier transform an MRA
$(\mathcal{D},\mathcal{T})$ wavelet set.
\end{theorem}

\begin{theorem}\cite{RS}
If $\psi $ is a classical wavelet and the set
$K=\mathrm{supp}(\hat{\psi})$ satisfies

\begin{enumerate}
\item $\sum_{k\in \mathbb{Z}}\chi _{K}(\xi +k)\leq 2$ a.e.;

\item $\sum_{k\in \mathbb{Z}}\chi _{K}(2^{j}\xi )\leq 2$
a.e.
\end{enumerate}

Then $K$ contains a wavelet set.
\end{theorem}

Qing Gu has an unpublished example which shows that the
techniques
in \cite{RS} do not extend to the case that $\sum_{k\in
\mathbb{Z}}\chi_{K}(\xi +k)\leq 3$ a.e. and $\sum_{j\in
\mathbb{Z}}\chi_{K}(2^{j}\xi )\leq 3$ a.e.

Tilings and spectral sets are related by the \textit{Fuglede
conjecture} \cite{F74}

\begin{conjecture}[Fuglede]
A measurable  set $\Omega $, with positive and finite
measure is a
spectral set if and only if $\Omega $ is an additive $\cT$
tile
for some set $\cT$.
\end{conjecture}

The conjecture, in general,  still remains unsolved, even if several
partial
results have been obtained \cite{IKT99,LW97,JP98,JP99,W02}.
In June 2003 it was shown by Tao,
\cite{TT03} that the conjecture in false in dimension
$5$ and higher.
We
will not discuss those articles,   but concentrate on the important paper
\cite{W02} by
Wang, which also made the first serious  attempt at studying
$(\mathcal{D},\mathcal{T})$
wavelet sets when $\mathcal{D}$ is not even a subgroup of
$\GL(n,\mathbb{R})$,
and $\mathcal{T}$ is not a lattice. We need two more
definitions before we
state some of Wang's results. Let $a\in \GL(n,\mathbb{R})$.
A set $\mathcal{D}\subseteq \GL(n,\mathbb{R})$ is said to be $a$ invariant if
$\mathcal{D}a=%
\mathcal{D}$. The multiplicative tiling set $\mathcal{D}$
said to
satisfy the \textit{interior condition} if there exists a
multiplicative $\mathcal{D}$-tile $\Omega $
such that $\Omega ^{o}\not=\emptyset $. Similarly the
spectrum $%
\mathcal{T}\subset \mathbb{R}^{n}$ satisfies the
\textit{interior
condition}
if there exists a measurable set $\Omega \subset
\mathbb{R}^{n}$ such that $%
\Omega ^{o}\not=\emptyset $ and $(\Omega ,\mathcal{T})$ is a
spectral pair. With these definitions we can state two of
Wang's
main results:

\begin{theorem}[Wang,\cite{W02}] \label{wset}
Let $\mathcal{D}\subset \GL(n,\mathbb{R})$ and
$\mathcal{T}\subset \mathbb{R}%
^{n}$. Let $\Omega \subset \mathbb{R}^{n}$ be measurable,
with
positive and finite measure. If $\left\{ d^{T}(\Omega )\mid
d\in
\mathcal{D}\right\} $ is a tiling of $\mathbb{R}^{n}$ and
$(\Omega ,\mathcal{T})$ is a spectral pair,
then $\Omega $ is a $(\mathcal{D},\mathcal{T})$-wavelet set.
Conversely, if $%
\Omega $ is a $(\mathcal{D},\mathcal{T})$-wavelet set and
$0\in \mathcal{T}$,
then $\Omega $ is a multiplicative $\mathcal{D}^T$-tile
and $(\Omega ,%
\mathcal{T})$ is a spectral pair.
\end{theorem}

\begin{theorem}[Wang,\cite{W02}]
\label{wang21} Let $\mathcal{D}\subset
\mathrm{GL}(n,\mathbb{R})$ such that $%
\mathcal{D}^{T}:=\{d^{T}\mid d\in \mathcal{D}\}$ is a bounded
multiplicative tiling set,
and let $\mathcal{T}\subset \mathbb{R}^{n}$ be a spectrum,
with both $%
\mathcal{D}^{T}$ and $\mathcal{T}$ satisfying the interior
condition. Suppose that $\mathcal{D}^{T}$ is $a$-invariant
for
some expanding matrix $a$
and $\mathcal{T}-\mathcal{T}\subset \mathcal{L}$ for some
lattice $\mathcal{L%
}$ of $\mathbb{R}^{n}$. Then, there exists a wavelet set
$\Omega $
with respect to $\mathcal{D}$ and $\mathcal{T}$.
\end{theorem}

In his paper, Wang states ``The assumption
that
$\cD^T$ ... have the interior condition is most likely
unnecessary.  All known examples of multiplicative tiling
sets
admit a tile having nonempty interior.''
In this section, we will
in fact show that the assumption that $\mathcal{D}^{T}$
satisfies
the interior condition is indeed unnecessary, but not by
proving
that every multiplicative tiling set admits a tile having
nonempty
interior. Instead, we will use a Lebesgue density argument
as in
\cite{DLSII, DDGH}. Moreover, the assumption of
multiplicative
tiling sets having prototiles that are bounded and bounded
away
from the origin is not a ``wavelet" assumption, but rather
it is
motivated from the point of view of tiling questions and the
relation between translation and dilation tilings of the
line.
From the point of view of wavelets, by Theorem \ref{wset},
one
does not always wish to restrict to bounded multiplicative
tiling
sets.  There are, however, some benefits of obtaining
wavelet sets
that are bounded and bounded away from the origin -
especially if
they also satisfy some additional properties.   For example,
if
the sets are the finite union of intervals, one can use
these
wavelets to show that theorems about the poor decay of
wavelets in
$L^{2}(\mathbb{R}^{n})$ for \textquotedblleft bad" dilations
are
optimal. Along these lines, Bownik \cite{Bo} showed that if
$a$ is
irrational and $\psi_1,\ldots,\psi_L$ is an $(A, \Z)$
multiwavelet, then there is an $i$ such that for each
$\delta >
0$, $\limsup_{\abs{x}\to \infty} \abs {\psi_i}{\abs {x}
}^{1+\delta} = \infty$.  He also showed that this result is
sharp
by finding wavelet sets for each of these dilations that are
the
union of at most three intervals.  Another possibility is to
use
wavelet sets that are the finite union of intervals (and
satisfying several extra conditions) as a start point for
the
smoothing techniques in \cite{DL, HW} However, these two
advantages come from having wavelet sets that are not only
bounded
and bounded away from the origin, but also the finite union
of
intervals. In the construction considered in \cite{W02}, it
is not
clear at all whether the end wavelet sets can be chosen to
be the
finite union of nice sets. In fact, the construction used of
Benedetto and Leon was used originally exactly to construct
fractal-like wavelet sets.

Since the general question of existence of wavelet sets is
phrased
not in terms of sets bounded and bounded away from the
origin, but
arbitrary measurable sets, we will also show that the
assumption
that there exist a multiplicative tiling set that is bounded
and
bounded away from the origin is unnecessary. This will be
done by
showing that whenever there is a set that tiles
$\mathbb{R}^n$ by
$\mathcal{D}$ dilations, where $\mathcal{D}$ is invariant
under an
expansive matrix, then there exists a bounded multiplicative
tiling set for $\mathcal{D}$.

We begin with some easy observations that were also in
\cite{W02}.
We say that sets $U$ and $V$ in $\mathbb{R}^n$ are
$a$-dilation
equivalent if there is a partition $\{U_k: k\in
\mathbb{Z}\}$ of
$U$ such that $\{a^k U_k \mid  k\in \mathbb{Z}\}$ is a partition
of
$V$.

\begin{lemma}
\label{invariant} Let $\mathcal{D} \subset \GL(n,
\mathbb{R})$
invariant under the invertible matrix $a$. If $\{d\Omega\mid d\in
\mathcal{D}\}$ is a tiling of $\mathbb{R}^n$ and $\Omega_0$
is
$a$-dilation equivalent to $\Omega$, then $\{d\Omega_0\mid  d\in
\mathcal{D}\}$ is a tiling of $\mathbb{R}^n$.
\end{lemma}

\begin{proof}
Let $S_k$ be a partition of $\Omega$ such that $\Omega_0 =
\bigcup_{k\in \mathbb{Z}} a^k (S_k)$. Then,
\begin{eqnarray*}
\bigcup_{d\in\mathcal{D}} d\Omega_0 &=& \bigcup_{D\in
\mathcal{D}}
\bigcup_{j\in \mathbb{Z}} da^j S_j\cr &=&
\bigcup_{j\in\mathbb{Z}}
\bigcup_{d\in \cD} da^j S_j\cr %
&=&\bigcup_{j\in \mathbb{Z}} \bigcup_{d\in \mathcal{D}} d
S_j =
\mathbb{R}^n.
\end{eqnarray*}
Similarly, one can show that $d_1 \Omega_0 \cap d_2\Omega_0$
has
measure 0 for all $d_1\not= d_2$ in $\mathcal{D}$.
\end{proof}

\begin{lemma} \label{equiv} Let $A$ be an expansive matrix
and
$\Omega_0, \Omega_1$ be such that $\abs{A^j\Omega_i \cap
A^k\Omega_i} = \delta_{j,k}$ for $i = 1,2$.  Then,
$\Omega_0$ is
$A$ equivalent to $\Omega_1$ if and only if $\cup_{j\in \Z}
A^j
\Omega_0 = \cup_{j\in \Z} A^j \Omega_1$ a.e.
\end{lemma}

\begin{proposition}
\label{bounded} Let $\mathcal{D} \subset \GL(n, \mathbb{R})$
be
invariant under the expansive matrix $a$. If there is a set
$\Omega\subset \mathbb{R}^n$ such that $\{d\Omega\mid  d\in
\mathcal{D}\}$ tiles $\mathbb{R}^n$, then there is a set
$\Omega_0$ bounded and bounded away from the origin such
that
$\{d\Omega_0\mid  d\in \mathcal{D}\}$ tiles $\mathbb{R}^n$. In
particular, $%
\mathcal{D}$ is a bounded multiplicative tiling set.
\end{proposition}

\begin{proof}
It is widely known that $a$ is expansive if and only if
there is
an ellipsoid $\mathcal{E}$ such that $\overline{\mathcal{E}}
\subset a \mathcal{E}^\circ$. In this case, it is easy to
check
that $\Omega_1 = a\mathcal{E} \setminus \mathcal{E}$ is a
bounded
multiplicative tiling set for $\{a^j\mid  j\in \mathbb{Z}\}$;
that is,
$\Omega_1$ is bounded and bounded away from the origin, and
$\{a^j\Omega_1\mid  j\in \mathbb{Z}\}$ tiles $\mathbb{R}^n$. Let
$S_j
= a^j(\Omega_1) \cap \Omega$, and $\Omega_0 = \bigcup_{j\in
\mathbb{Z}} a^{-j} S_j$. It is clear that $\Omega_0\subset
\Omega_1$, so it is bounded and bounded away from the
origin.
Moreover, since $\{a^j(\Omega_1) \mid  j\in \mathbb{Z}\}$ is a
tiling
of $\mathbb{R}^n$, it follows that $\{S_j \mid  j\in
\mathbb{Z}\}$ is
a partition of $\Omega$; hence, $\Omega_0$ is $a$-dilation
equivalent to $\Omega$. Therefore, by lemma \ref{bounded},
$\Omega_0$ is a multiplicative tiling set.
\end{proof}

Next, we turn to showing that the assumption of a
multiplicative
tile with non-empty interior is unnecessary. We have
(combining
Lemma 2 and Theorem 1 of \ref{ddgh}):

\begin{theorem}[\cite{DDGH}]
\label{ddgh} Let $M$ be a measurable subset of
$\mathbb{R}^n$ with
positive measure satisfying $a M = M$ for some expansive
matrix
$a$.
Then, there exists a set $E \subset M$ such that $\{E + k \mid
k\in \mathbb{Z}%
^n\}$ tiles $\mathbb{R}^n$ and $\{a^j E \mid  j\in \mathbb{Z}\}$
tiles
$M$.
\end{theorem}

Suppose that we are considering classes of $(\mathcal{D},
\mathcal{T})$ wavelets, where $\mathcal{D} = \{a^j\mid
j\in\mathbb{Z}\}$ and $\mathcal{T}$ is a lattice. It is a
general
principle that one can either assume that $a$ is in (real)
Jordan
form, in which case one must deal with arbitrary lattices,
or one
can assume that the lattice $\mathcal{T} = \mathbb{Z}^n$, in
which
case one needs to consider all matrices of the form
$bab^{-1}$. In
particular, if one is working with expansive matrices, it is
almost always permissible to restrict attention to
translations by
$\mathbb{Z}^n$. While this is clear to experts in the field,
it is
likely that researchers new to this field are not aware that
the
above theorem is really a theorem about arbitrary lattices.

Indeed, let $M$ be a measurable subset of $\mathbb{R}^n$
with
positive measure satisfying $a M = M$, for some expansive
matrix
$a$. Let $\cL$ be a full rank lattice in $\mathbb{R}^n$.
Then,
there is an invertible matrix $b$ such that $b\cL =
\mathbb{Z}^n$.
The set $b M$ is $bab^{-1}$ invariant, and $bab^{-1}$ is an
expansive matrix, so there is a set $F$ such that $\{F+ k \mid
k\in
\mathbb{Z}^n\}$ tiles $\mathbb{R}^n$ and $\{ba^jb^{-1}F \mid
j\in
\mathbb{Z}\}$ tiles $b M$. We claim that for $E = b^{-1}F$,
$\{E +
k \mid  k\in \cL\}$ tiles $\mathbb{R}^n$ and $\{a^jE \mid  j\in
\mathbb{Z}\}$ tiles $M$. Indeed,
\begin{eqnarray*}
\bigcup_{k\in \cL} E + k &=&
\bigcup_{k\in \cL} b^{-1} F + k\cr
&=& \bigcup_{k\in \cL} b^{-1}(F + b k)\cr
&=& \bigcup_{k\in \mathbb{Z}^n} b^{-1}(F + k)\cr
& =& b^{-1}\bigl(\bigcup_{k\in \mathbb{Z}^n} (F + k)\bigr)\\
&= &\mathbb{R}^n\, .
\end{eqnarray*}
One can similarly show the disjointness of translates by
$\cL$. To
see that $\{a^jE \mid  j\in \mathbb{Z}\}$ tiles $M$, note that
\begin{eqnarray*}
\bigcup_{j\in \mathbb{Z}} a^j E &=& \bigcup_{j\in
\mathbb{Z}} a^j b^{-1}F\cr %
&=& b^{-1}\bigl(\bigcup_{j\in \mathbb{Z}} ba^jb^{-1}
F\bigr)\cr
&=& b^{-1} b M = M.
\end{eqnarray*}
Again, disjointness of the dilates is immediate. Thus, we
have
proven the following theorem, that seems to be well known:

\begin{theorem}[\cite{DDGH}]
\label{ddgh2} Let $M$ be a measurable subset of
$\mathbb{R}^n$
with positive measure satisfying $a M = M$ for some
expansive
matrix $a$, and let $\mathcal{T}$ be a full rank lattice in
$\mathbb{R}^n$. Then, there
exists a set $E \subset M$ such that $\{E + k \mid  k\in
\mathcal{T}\}$ tiles $%
\mathbb{R}^n$ and $\{a^jE \mid  j\in \mathbb{Z}\}$ tiles $M$.
\end{theorem}

Theorem \ref{ddgh2} can be used to give an easy proof of
Theorem
\ref{wang21} removing three of the assumptions, but adding
the
assumption that the translation set is a lattice.

\begin{theorem}
\label{improvement} Let $\mathcal{D} \subset \GL(n,
\mathbb{R})$
be such that there exists a measureable set $\Omega$ such
that
$\{d^T \Omega\mid d\in \mathcal{D}\}$ is a tiling of
$\mathbb{R}^n$. Suppose also that $\mathcal{D}^T$ is $a$
invariant
for some expansive matrix $a$. Let $\mathcal{T} \subset
\mathbb{R}^n$ be a full-rank lattice. Then, there exists a
($\mathcal{D}, \mathcal{T})$ wavelet set $E$.
\end{theorem}

\begin{proof}
By assumption, there exists a set $\Omega$ such that
$\mathcal{D}^T(\Omega)$ is a tiling of $\mathbb{R}^n$.
Consider
the set $M = \bigcup_{j\in \mathbb{Z}} a^j \Omega$. The set
$M$ is
clearly $a$ invariant, and $\{a^j(\Omega)\mid j\in \Z\}$ is
a
measurable partition of $M$ so by \ref{ddgh2}, there exists
a set
$E$ such that $\{a^j(E) \mid j\in \mathbb{Z}\}$ tiles $M$
and $\{E
+ k \mid  k\in \mathcal{L}^*\}$ tiles $\mathbb{R}^n$. By Lemmas
\ref{invariant} and \ref{equiv}, since $E$ is $a$-equivalent
to
$\Omega$, $\{d^TE\mid d\in \mathcal{D}\}$ tiles
$\mathbb{R}^n$.
That is, $E$ is a $(\mathcal{D}, \mathcal{T})$ wavelet set,
as
desired.
\end{proof}

We have exhibited above the essential nature of the argument
in
\cite{W02}. That is, what is desired is a general criterion
for
the following question:
\begin{question}
Given an expansive matrix $a$, a lattice $\mathcal{L}$ and
two
sets $\Omega_1 $ and $\Omega_2$, when does there exist a set
$\Omega$ that is $a$-equivalent to $\Omega_1$ and
$\mathcal{L}$
equivalent to $\Omega_2$?
\end{question}
In the above case, we were forced to restrict to the case
that
$\Omega_2$ is a fundamental region for the lattice
$\mathcal{L}$,
since that is what was shown in \cite{DDGH}. As a final
generalization in this section, we show that what is really
necessary is that $\Omega_2$ contain a neighborhood of the
origin.
The reader should compare the theorem below with the
statement and
proofs of the theorems in \cite{DLS} and \cite{DLSII}.

\begin{theorem}\label{dlsimprovement}  Let $a$ be an
expansive matrix and
$\Omega_1\subset \brn$ a set of positive measure such that
$|a^j\Omega_1 \cap a^k \Omega_1| = 0$ whenever $j\not= k$.
Let $M
= \bigcup_{j\in \Z} a^j \Omega_1$.  Let $\cL\subset \brn$ be
a
full rank lattice and $\Omega_2\subset \brn$ such that
$|\Omega_2
+ k_1 \cap \Omega_2 + k_2| = 0$ for $k_1\not=k_2 \in \cL$
and such
that there exists $\epsilon > 0$ such that $M \cap
B_\epsilon(0)
\subset \Omega_2 \cap B_\epsilon(0) $.  Then, there exists a
set
$\Omega$ such that $\Omega$ is $a$ equivalent to $\Omega_1$
and
$\cL$ equivalent to $\Omega_2$.
\end{theorem}

Before proving Theorem \ref{dlsimprovement}, we state and
prove
its main corollary, which is Theorem 2.1 of \cite{W02} with
all
but one technical assumption removed.

\begin{corollary}
\label{bestwang} Let $\mathcal{D} \subset \GL(n,
\mathbb{R})$ be
such that $\cD^T$ is a multiplicatvie tiling set. Let
$\mathcal{T}$ be a spectrum with interior such that there
exists a
full rank lattice such that $\mathcal{T} - \mathcal{T}
\subset
\mathcal{L}$. Then, if $\mathcal{D}^T$ is $a$-invariant for
some
expansive matrix $a$, there exists a $(\mathcal{D},
\mathcal{T})$
wavelet set.
\end{corollary}

\begin{proof}
Since translations of spectral sets are again spectral sets,
we
may assume without loss of generality that $\Omega_2$
contains $0$
as an interior point. By Lemma 3.1 of \cite{W02}, $(\Omega_2 +
k_1) \cap (\Omega_2 + k_2)$ has measure $0$ whenever
$k_1\not=k_2\in \mathcal{L}^*$. So, by \ref{dlsimprovement},
there
is a set $\Omega$ that is $a$ equivalent to $\Omega_1$ and
$\mathcal{L}^*$ equivalent to $\Omega_2$. By Lemma
\ref{invariant}, $\{d^T\Omega\mid  d\in \mathcal{D}\}$ tiles
$\mathbb{R}^n$, and by Lemma 3.2 in \cite{W02}, $(\Omega,
\mathcal{T})$ is a spectral set. Therefore, $\Omega$ is a
$(\mathcal{D}, \mathcal{T})$ wavelet set.
\end{proof}

We turn now to proving Theorem \ref{dlsimprovement}.  We
begin by
noting that arguing as in the proof of Theorem \ref{ddgh2},
one
can restrict to the case that $\cL = \Z^n$.  Next, we need
to
extract the following lemma from the proof of Corollary 1 in
\cite{DDGH}, then we will follow very closely the proof in
\cite{DLSII}.

\begin{lemma}\label{ddghlemma} Let $a$ be an expansive
matrix in $GL(n, \R)$. Let $F_0$ be a set of
positive measure such that $|a^jF_0 \cap a^k F_0| = 0$
whenever
$j\not= k$.   Let $E = [-1/2,1/2]^n$.  Then, for every
$\epsilon >
0$, there exists $k_0\in \Z$ and $\ell_0\in \Z^n$ such that
$|a^{k_0} F_0 \cap (E + \ell_0)| > 1-\epsilon.$
\end{lemma}

The proof of Lemma \ref{ddghlemma} is a clever use of a
Lebesgue
density argument, which we will not repeat here.

%
%

\begin{proof}[Proof of Theorem \ref{dlsimprovement}]  First,
note that as in the case of Theorem
\ref{ddgh2}, Lemma \ref{ddghlemma} is really a lemma about
arbitrary full-rank lattices $\cL$.   Moreover, one can
replace $E
= [-1/2, 1/2]^n$ by any subset $E$ of a fundamental region
of
$\cL$ to get the following formally stronger lemma.

\begin{lemma}\label{ddghstronglemma} Let $a$ be an expansive
matrix in $GL(n, \R)$. Let $F_0$ be a set of
positive measure such that $|a^jF_0 \cap a^k F_0| = 0$
whenever
$j\not= k$.   Let $\cL\subset \brn$ be a full-rank lattice
with
fundamental region $\Omega$.  Then, for every set $E\subset
\Omega$ and every $\epsilon > 0$, there exists $k_0\in \Z$
and
$\ell_0\in \cL$ such that $|a^{k_0} F_0 \cap (E + \ell_0)|
\ge
(1-\epsilon)|E|.$
\end{lemma}

Turning to the proof of \ref{dlsimprovement}, note that by
Theorem
\ref{bounded}, we may assume without loss of generality that
$\Omega_1$ is bounded and bounded away from the origin.   We
may
also assume that $\Omega_2$ is contained in a convex,
centrally
symmetric fundamental region of $\cL$.  The rest of the
proof
follows very closely the proof of Theorem 3.7 in
\cite{DLSII},
with Lemma \ref{ddghstronglemma} playing the role of
Proposition
3.5 in \cite{DLSII}.   We will construct a family $\{G_{ij}\mid
i\in
{\mathbb {N}}, j\in \{1,2\}\}$ of measurable sets whose
$a$-dilates form a measurable partition of $\Omega_1$ and
whose
translates by vectors in $\cL$ form a measurable partition
of
$\Omega_2$.  Then
\begin{equation}
\Omega:=\bigcup G_{ij}
\end{equation}
is the set desired in Theorem \ref{dlsimprovement}.  Since
the
steps are so similar to \cite{DLSII}, we will give the first
step
of the inductive definition, and the properties needed for
induction.  Details are the same as in \cite{DLSII}.

Let $\{\alpha_i\}$ and $\{\beta_i\}$ be sequences of
positive
constants decreasing to 0 and such that $\alpha_1 <
\epsilon$
chosen so that $B_{\alpha_1}(0)  \cap M \subset
B_{\alpha_1}(0)
\cap \Omega_2$.  Let $\tilde E_{11} = \Omega_2\setminus
B_{\alpha_1}(0)$.  Then, $|(\tilde E_{11})| > 0.$  Let
$\tilde
F_{11}$ be a measurable subset of $\Omega_1$ with measure
strictly
less than $|\Omega_1|$.  By Lemma \ref{ddghstronglemma},
there
exists $k_1\in \N$ and $\ell_1\in \cL$ such that
\begin{equation}
|a^{k_1} \tilde F_{11} \cap (\tilde E_{11} - \ell_1)| \ge
\frac{1}{2} |\tilde E_{11}|
\end{equation}
Let $G_{11} := a^{k_1} \tilde F_{11} \cap (\tilde E_{11} -
\ell_1),$ let $E_{11} := G_{11} + \ell_1$, and let
\begin{equation}
F_{11} := \tilde F_{11} \cap a^{-k_1} (E_{11} - \ell_1) =
a^{-k_1}
G_{11}
\end{equation}
Then $F_{11} \subset \tilde F_{11} \subset \Omega_1$ and
$|\Omega_1\setminus \tilde F_{11}| \ge |\Omega_1 \setminus
F_{11}|
>  0$.  Also, $E_{11} \subset \tilde E_{11},$ and
\begin{equation}
|E_{11}| = |G_{11}| \ge \frac{1}{2} |\tilde E_{11}|
\end{equation}
Also, $G_{11} = a^{k_1} F_{11}$.  Now choose $F_{12}\subset
\Omega_1$, disjoint from $F_{11}$, such that
$\Omega_1\setminus(F_{11}\bigcup F_{12})$ has positive
measure
less than $\beta_1$.  Choose $m_1$ such that $a^{-m_1}
F_{12}$ is
contained in $N_1 = B_{\alpha_1/2}(0)$ and is disjoint from
$G_{11}$.  (This is possible since $G_{11}$ is bounded away
from
0.)  Set
\begin{equation}
G_{12}:=E_{12}:=a^{-m_1} F_{12}.
\end{equation}
The first step is complete.

Proceed inductively, obtaining disjoint families of positive
measure $\{E_{ij}\}$ in $\Omega_2$, $\{F_{ij}\}$ in
$\Omega_1$ and
$\{G_{ij}\}$ such that for $i=1,2,\ldots$ and $j=1,2$ we
have

\begin{enumerate}
\item $G_{i1} + \ell_i = E_{i1}$;
\item $G_{i2} = E_{i2}$;
\item $a^{-k_1} G_{i1} = F_{i1}$;
\item $a^{m_i} G_{i2} = F_{i2}$;
\item $|\Omega_1\setminus (F_{11} \bigcup F_{12}
\bigcup\cdots \bigcup F_{i1} \bigcup F_{i2})| < \beta_i,$
and
\item $|E_{i1}| \ge \frac 12 |\Omega_2\setminus (E_{11}
\bigcup E_{12} \bigcup\cdots \bigcup E_{i-1,1} \bigcup
E_{i-1,2})) - \frac 12 |N_i|,$ where $N_{i}$ is
a ball centered at 0 with radius less than $\alpha_i$.
\end{enumerate}

Since $\beta_i\to 0$, item 5 implies that $F\setminus
(\bigcup
F_{ij})$ is a null set, and since $\alpha_i\to 0$, item 6
implies
that $(E\setminus (\bigcup E_{ij}))$ is a null set.  Let
\begin{equation}
F = \bigcup\{G_{ij} \mid  i=1,2,\ldots,j = 1,2\}
\end{equation}
then, $G$ is congruent to $\Omega_2$ by items 1 and 2, and
the $a$
dilates of $G$ form a partition of $M$, as desired.
\end{proof}

For sets $\mathcal{D}\subset \GL(n, \mathbb{R})$ which are
invariant under an expansive matrix, Theorem \ref{bestwang}
is
nice in that it reduces the question of existence of wavelet
sets
to the the question of existence of tiling sets for
dilations and
translations separately. It is still in some sense
unsatisfactory,
because it relies on the existence of objects external to
the sets
$(\mathcal{D}, \mathcal{T})$ under consideration. From the
point
of view of characterizing sets $(\mathcal{D}, \mathcal{T})$
for
which wavelet sets exist, something more is needed.  We will
present in section 4 some progress on this question when
$\cD$ is
a countable subgroup of $GL(n, \R)$.

%
%

\section{Admissible groups and frames}

\noindent We will now turn to the applications of those
results to
frames and wavelets in $\R^n$. But first we recall some
results
about the continuous wavelet transform.

Recall that translations and dilations on the real line form
the
so-called $(ax+b)$-group. Assume in general that we have a
group
$G$ acting on a topological space. Assume that $\mu$ is a
Radon
measure on $X$ such that $\mu$ is \textit{quasi-invariant},
i.e,
there exists a measureable function $j:G\times X\to \R^+$
such
that
$$\int_X f(g\cdot x)\, d\mu (x)=\int_X j(g,x)f(x)\, d\mu
(x)$$
for all $f\in L^1(X)$. Then we can define a unitary
representation
of $G$ on $L^2(X)$ by
$$[\pi (g)f](x)=j(g^{-1},x)^{-1/2}f(g^{-1}x)\, .$$
For a fixed $\psi \in L^2(G)$ define the transform
$W_\psi:L^2(X)\to C(G)$ by
$$W_\psi (f)(g):=(f,\pi (g)\psi)$$
and notice that $W_\psi$ intertwines the representation
$\pi$ and
the \textit{left regular representation}, i.e.,
\begin{eqnarray*}
W_\psi (\pi(g)f)(x)&=&(\pi (g)f,\pi (x)\psi)\\
&=&(f,\pi (g^{-1}x)\psi)\\
&=&W_\psi (f)(g^{-1}x)\, .
\end{eqnarray*}

For the $(ax+b)$-group this becomes

\begin{equation}  \label{contiwave}
W_\psi (f)(a,b)= (f\mid \pi
(a,b)\psi)=|a|^{-1/2}\int_{\mathbb{R}}f(x)%
\overline{\psi ((x-b)/a )}\, dx\, .
\end{equation}
Here $T_b:L^2(\mathbb{R})\to L^2(\mathbb{R})$ stands for the
unitary
isomorphism corresponding to translation $T_bf(x)=f(x-b)$
and $D_a:L^2(%
\mathbb{R})\to L^2(\mathbb{R})$ is the unitary map
corresponding
to dilation $D_af(x)=|a|^{-1/2}f(x/a)$, $a\not=0$.

The discrete wavelet transform is obtained by sampling the
wavelet
transform, given by a suitable wavelet $\psi$, of a function
$f$
at points gotten by replacing the full $(ax+b)$-group by a
discrete subset generated by translation by integers and
dilations
of the form $a=2^n$:
\begin{eqnarray*}
W_\psi^d(f)(2^{-n},-2^{-n}m) &=&(f\mid \pi
((2^n,m)^{-1})\psi ) \\
&=&2^{n/2}\int_{\mathbb{R}}f(x)\overline{\psi (2^n x+m)}\,
dx\, .
\end{eqnarray*}
Hence, the corresponding frame is
\begin{equation}  \label{frame}
\{\pi ((2^n,m)^{-1})\psi\mid n,m\in \mathbb{Z}\}\, .
\end{equation}
The inverse refers here to the inverse in the
$(ax+b)$-group.

In the same way it is well known, that the short time
Fourier
transform, and several other well known integral transforms
have a
common explanation in this way in the language of
representation
theory.   This observation is the basis for the
generalization of
the continuous wavelet transform to higher dimensions and
more
general settings, and was already made by A.
Grossmann, J. Morlet, and T. Paul in 1985, see
\cite{GMP85,GMP86}. In \cite%
{GMP85} the connection to square integrable representations
and
the relation to the fundamental paper of M. Duflo and C. C.
Moore
\cite{DM76} was already pointed out. Several natural
questions
arise now, in particular to describe the image of the
transform
$W_\psi$ and how that depends on $\psi$. But we will not go
into
that here, but refer to
\cite{AAG2000,ALTW2002,BCMO95,BT96,DGM86,FO02,FG92,FG89,%
F96,F00,FM01,G98,%
G2001,GMP85,GMP86,HW89,LWWW2002,O02,WW01} for discussion.
Here, we
will concentrate on the connection to frames, wavelets and
wavelet
sets.

Denote by $\mathrm{Aff}(\mathbb{R}^{n})$ the group of
invertable
affine linear transformations on $\mathbb{R}^{n}$. Thus
$\mathrm{Aff}(%
\mathbb{R}^{n})$ consists of pairs $(x,h)$ with $h\in
\mathrm{GL}(n,\mathbb{R%
})$ and $x\in \mathbb{R}^{n}$. The action of $(x,h)\in
\mathrm{Aff}(\mathbb{R%
}^{n})$ on $\mathbb{R}^{n}$ is given by
\begin{equation*}
(x,h)(v)=h(v)+x\,.
\end{equation*}%
The product of group elements is the composition of maps.
Thus
\begin{equation*}
(x,a)(y,b)=(a(y)+x,ab)
\end{equation*}%
the identity element is $e=(0,\mathrm{id})$ and the inverse
of
$(x,a)\in \mathrm{Aff}(\mathbb{R}^{n})$ is given by
\begin{equation*}
(x,a)^{-1}=(-a^{-1}(x),a^{-1})\,.
\end{equation*}%
Thus $\mathrm{Aff}(\mathbb{R}^{n})$ is the semidirect
product of
the abelian
group $\mathbb{R}^{n}$ and the group
$\mathrm{GL}(n,\mathbb{R})$; $\mathrm{%
Aff}(\mathbb{R}^{n})=\mathbb{R}^{n}\times
_{s}\mathrm{GL}(n,\mathbb{R})$. Let $H\subset
\mathrm{GL}(n,\mathbb{R})$ be a closed subgroup. Define
\begin{equation*}
\mathbb{R}^{n}\times _{s}H:=\{(x,a)\in
\mathrm{Aff}(\mathbb{R}^{n})\mid a\in H,~x\in
\mathbb{R}^{n}\}\text{.}
\end{equation*}%
Then $H\times _{s}\mathbb{R}^{n}$ is a closed subgroup of
$\mathrm{Aff}(%
\mathbb{R}^{n})$.

Define a unitary representation of $\mathbb{R}^{n}\times
_{s}H$ on $L^{2}(%
\mathbb{R}^{n})$ by
\begin{equation}
\lbrack \pi (x,a)f](v):=|\det
(a)|^{-1/2}f((x,a)^{-1}(v))=|\det
(a)|^{-1/2}f(a^{-1}(v-x))\,.
\end{equation}%
Write $\psi _{x,a}$ for $\pi (x,a)\psi $. We will
also need
another action of $H$ on $\mathbb{R}^{n}$ by
$a\cdot \omega :=(a^{-1})^{T}(\omega )$.
We denote by $\hat{\pi}(x,a)$ the
unitary action on $L^{2}(\mathbb{R}^{n})$ given by
\begin{equation}
\hat{\pi}(x,a)f(v)=\sqrt{|\det (a)|}e^{-2\pi i(x\mid
v)}f(a^{-1}\cdot v)=%
\sqrt{|\det (a)|}e^{-2\pi i(x\mid v)}f(a^{T}(v))\,.
\label{rep-fourier}
\end{equation}%
\begin{remark} Some authors use the semidirect product
$H\times_s\R^n$ instead
of $R^n\times_s H$. Thus first the translation and then the
linear
map is applied, i.e., $(a,x)(v)=a(v+x)$. In this notation
the
product becomes $(a,x)(b,y)=(ab,ab(y)+a(x))$, the inverse of
$(a,x)$ is $(a,x)^{-1}=(a^{-1},-ax)$, and the wavelet
representation is
$$\widetilde{\pi}(a,x)f(v)=|\det a|^{-1/2}f(a^{-1}v-x)\, . $$
\end{remark}
The Fourier transform intertwines the representations $\pi $
and
$\hat{\pi}$ \cite{FO02}, Lemma 3.1, i.e.,
\begin{equation}
\widehat{\pi (x,a)f}(\omega )=\hat{\pi}(x,a)\hat{f}(\omega
)\,,\quad f\in L^{2}(\mathbb{R}^{n})\,.  \label{fourier}
\end{equation}%
Denote by $d\mu _{H}$ a left invariant measure on $H$. A
left
invariant measure on $G$ is then given by
$d\mu_{G}(x,a)=|\mathrm{\det }(a)|^{-1}d\mu _{H}(a)dx$. Let
$f,\psi \in
L^{2}(\mathbb{R}^{n})$. A simple calculation shows that
\begin{eqnarray}
\int_{G}|(f\mid \pi (x,a)\psi )|^{2}\,d\mu _{G}(x,a)
&=&\int_{\mathbb{R}%
^{n}}|\hat{f}(\omega )|^{2}\int_{H}|\hat{\psi}(a^{-1}\cdot
\omega
)|^{2}\,d\mu _{H}(a)d\omega  \label{eq-integra} \\
&=&\int_{\mathbb{R}^{n}}|\hat{f}(\omega
)|^{2}\int_{H}|\hat{\psi}%
(a^{T}(\omega ))|^{2}\,d\mu _{H}(a)d\omega
\end{eqnarray}%
There are several ways to read this. First let $M\subseteq
\R^n$
be measurable and invariant under the action $H\times
\R^n\ni
(h,v)\mapsto h\cdot v:=(h^{-1})^T(v)\in \R^n$. Then we
denote by
$L^2_M(\R^n)$ the space of function $f\in L^2(\R^n)$ such
that
$\hat{f}(\xi)=0$ for almost all $\xi \not\in M$. Notice that
$L^2_M(\R^n)$ is a closed invariant subspace, and that the
orthogonal projection onto $L^2_M(\R^n)$ is given by
$f\mapsto
(\hat{f}\chi_M)^\vee$. The first result is now:

\begin{theorem} Let $M\subseteq \R^n$ be measurable of
positive measure, and invariant under the action
$(a,v)\mapsto a\cdot v$. Then the wavelet transform
\begin{equation*}
W_{\psi }:f\mapsto (f\mid \pi (x,a)\psi
)_{L^{2}(\mathbb{R}^{n})}=|\det a|^{-1/2}\int
f(y)\overline{\psi
(a^{-1}(y-x))}~dy
\end{equation*}%
is a partial isometry $W_\psi :L^2_M(\R^n)\to L^{2}(G)$ if
and
only if
\begin{equation}
\Delta _{\psi }(\omega ):=\int_{H}|\hat{\psi}(a^{T}\omega
)|^{2}\,d\mu _{H}(a)=1  \label{eq-wavlet}
\end{equation}%
for almost all $\omega \in M$.
\end{theorem}

Following \cite{LWWW2002,WW01} we define
\begin{definition}[Laugesen, Weaver, Weiss, and Wilson]
Let $M\subseteq \R^n$, be measureable, invariant, and
$|M|>0$. Let
$\psi \in L^{2}(\mathbb{R}^{n})$ then $\psi $ is said to be
\textit{a (normalized) admissible $(H,M)$-wavelet} if for
almost
all $\omega \in
M$ we have%
\begin{equation*}
\int_{H}|\hat{\psi}(a^{T}\omega )|^{2}\,d\mu _{H}(a)=1~.
\end{equation*}%
We say that the pair $(H,M)$ is \textit{admissible} if a
$(H,M)$-admissible wavelet $\psi $ exists. If $M=\R^n$ then
we say
that $H$ is \textit{admissible} and that $\psi$ is a
(normalized)
\textit{wavelet} function.
\end{definition}

Assume that $\psi $ is a normalized admissible wavelet.
Then%
\begin{equation*}
G\ni (x,a)\mapsto F(x,a):=W_{\psi }f(x,a)\psi _{x,a}\in
L^{2}(\mathbb{R}^{n})
\end{equation*}%
is well defined and if $g\in L^{2}(\mathbb{R}^{n})$ then%
\begin{eqnarray*}
\int_{G}(F(x,a)\mid g)_{L^{2}(\mathbb{R}^{n})}~d\mu
_{G}(x,a)
&=&\int_{G}W_{\psi }f(x,a)\left( \int_{\mathbb{R}^{n}}\psi
_{x,a}(y)%
\overline{g(y)}~dy\right) ~d\mu _{G}(x,a) \\
&=&\int_{G}W_{\psi }f(x,a)\overline{W_{\psi }g(x,a)}~
d\mu_{G}(x,a) \\
&=&(W_{\psi }f\mid W_{\psi }g)_{L^{2}(G)} \\
&=&(f\mid g)_{L^{2}(\mathbb{R}^{n})}~.
\end{eqnarray*}%
Hence

\begin{lemma}
Assume that $\psi \in L^{2}(\mathbb{R}^{n})$ satisfies
$\int_{H}|\hat{\psi}%
(a^{T}\omega )|^{2}~d\mu _{H}(a)=1$. Then%
\begin{equation}
f=\int W_{\psi }f(x,a)\psi _{x,a}~d\mu _{G}(x,a)
\end{equation}%
as a weak integral for all $f\in L^{2}(\mathbb{R}^{n})$.
\end{lemma}

\begin{question}[Laugesen, Weaver, Weiss, and Wilson]
Give a characterization of admissible subgroups of
$\mathrm{GL}(n,\mathbb{R}%
) $.
\end{question}

It is easy to derive one necessary condition for
admissibility.
For $\omega \in \mathbb{R}^{n}$ let
\begin{equation}
H^{\omega }=\{h\in H\mid h\cdot \omega =\omega \}=\left\{
h\in
H\mid h^{T}(\omega )=\omega \right\}  \label{compstab}
\end{equation}%
be the stabilizer of $\omega $. Then admissibility implies
that
$H$ is compact for almost all $\omega \in \mathbb{R}^{n}$.
But
this condition is not sufficient and there is by now
\textit{no}
complete characterization of admissible group. The best
result up
to now is the following due to Laugesen, Weaver, Weiss, and
Wilson
\cite{LWWW2002}:

\begin{theorem}[Laugesen, Weaver, Weiss, and Wilson,
\cite{LWWW2002}]\label{t-LWWW}
Let $H$ be a closed subgroup of $\mathrm{GL}(n,\mathbb{R})$.
For
$\omega \in \mathbb{R}^{n}$ and $\epsilon >0$ let
\begin{equation*}
H_{\epsilon }^{\omega }:=\{h\in H\mid \| a\cdot \omega
-\omega \|
\leq \epsilon \}
\end{equation*}%
be the $\epsilon $-stabilizer of $\omega $. If either

\begin{enumerate}
\item  $G=\R^n\times_s H$ is non-unimodular, or

\item $\{\omega \in \mathbb{R}^n\mid H^\omega $ is
non-compact $\}$ has
positive Legesgue measure
\end{enumerate}
holds, then $H$ is not admissible. If both (1) and
\begin{enumerate}
\item[(3)] $\{\omega \in \mathbb{R}^n\mid H^\omega $ is
non-compact for all $
\epsilon>0\}$ has positive Lebesgue measure fail,
\end{enumerate}
then $H$ is admissible.
\end{theorem}

\section{A special class of groups}
\noindent In \cite{FO02} and \cite{O02} a special class of
groups
with finitely many open orbits were discussed. Those were
related
to the so-called \textit{prehomogeneous vector spaces of
parabolic
type} \cite{BR}. We start with a simple lemma:

\begin{lemma}
Let $H\subset \mathrm{GL}(n,\R)$ be a closed subgroup,
Assume that
$M\subseteq \R^n$ is up to set of measure zero a union of
finitely many open orbits $U_{1},\ldots ,U_{k}\subset
\mathbb{R}%
^{n}$ under the action $(a,\omega )\mapsto a\cdot \omega
=(a^{-1})^{T}(\omega )$. Assume furthermore that $H^{\omega
}$ is compact for $%
\omega \in U_{j}$, $j=1,\ldots ,k$. Then the pair $(H,M)$ is
admissible.
\end{lemma}

\begin{proof}
Fix $\omega _{j}\in U_{j}$. For $j=1,\ldots ,k$ let
$g_{j}\in
C_{c}(U_{j})$,
$g_{j}\geq 0$, $g\not=0$. Then the function%
\begin{equation*}
H\ni a\mapsto g_{j}(a^{T}(\omega _{j}))\in \mathbb{C}
\end{equation*}%
has compact support and $\int_{H}g_{j}(a^{T}
(\omega_{j}))~d\mu_{H}(a)>0$. Let $\omega \in U_{j}$. Choose $h\in H$ such
that
$\omega =h^{T}(\omega_{j})$.
This is possible because $U_{j}$ is homogeneous
under $H$. Then%
\begin{eqnarray*}
\int_{H}g_{j}(a^{T}(\omega ))~d\mu _{H}(a)
&=&\int_{H}g_{j}(a^{T}h^{T}(\omega _{j}))~d\mu _{H}(a) \\
&=&\int_{H}g_{j}((ha)^{T}(\omega _{j}))~d\mu _{H}(a) \\
&=&\int_{H}g_{j}(a^{T}(\omega _{j}))~d\mu _{H}(a)~.
\end{eqnarray*}%
Hence $\Delta _{j}=\int_{H}g_{j}(a^{T}(\omega ))~d\mu_{H}(a)>0$
is
independent of $\omega \in U_{j}$. Define $\varphi
:\mathbb{R}^{n}\rightarrow \mathbb{C}$ by%
\begin{equation*}
\varphi (\omega ):=\left\{
\begin{array}{cc}
g_{j}(\omega )/\Delta _{j} & \text{if }\omega \in U_{j} \\
0 & \text{if }\omega \not\in \bigcup_{j=1}^{k}U_{j}%
\end{array}%
\right.  ~.
\end{equation*}%
Then $\varphi \in C_{c}(\mathbb{R}^{n})$, so in particular
$\varphi \in L^{2}(\mathbb{R}^{n})$. Define
$\psi :=\varphi^{\vee
}$. Then $\psi $
satisfies the admissibility condition (\ref{eq-wavlet}) and
it follows that $H$ is admissible.
\end{proof}

\begin{question}
\label{q-openorbits}Classify all the admissible group with
finitely many open orbits of full measure.
\end{question}

\begin{lemma}\label{le-FH}
Let $H\subset \mathrm{GL}(n,\mathbb{R})$ be a closed
subgroup such
that $H$ can be written as $H=ANR=NAR=RAN$ with $R$
compact, $A$ simply connected abelian, and such that the
map%
\begin{equation*}
N\times A\times R\ni (n,a,r)\mapsto nar\in H
\end{equation*}%
is a diffeomorphism. Assume furthermore that $R$ and $A$
commute,
and that $R$ and $A$ normalize $N$. Finally assume that
there
exists a co-compact discrete subgroup $\Gamma _{N}\subset
N$. Let
$\Gamma _{A}\subset A$ be a
co-compact subgroup in $A$. Choose bounded measureable
subsets $\mathbb{F}%
_{A}\subset A$, and $\mathbb{F}_{N}\subset N$ such that
$N=\Gamma _{N}%
\mathbb{F}_{N}$, and $A=\Gamma _{A}\mathbb{F}$ and such that
the
union is
disjoint. Let $\Gamma =\Gamma _{A}\Gamma _{N}$ and
$\mathbb{F}_{H}=\mathbb{F}%
_{N}\mathbb{F}_{A}R\subset H\,$.  Then we have
\begin{equation*}
H=\bigcup_{\gamma \in \Gamma }\gamma \mathbb{F}_{H}
\end{equation*}%
and the union is disjoint. Furthermore we can choose
$\mathbb{F}_{H}$ such that $\mathbb{F}_{H}^{o}\not=
\emptyset$.
\end{lemma}

\begin{proof}
We have
\begin{eqnarray*}
\bigcup_{\gamma }\gamma \mathbb{F}_{N}\mathbb{F}_{A}R
&=&\Gamma_{A}\Gamma_{N}(\mathbb{F}_{N})\mathbb{F}_{A}R \\
&=&\Gamma _{A}N\mathbb{F}_{A}R\qquad
\mathrm{because}\,\Gamma _{N}\mathbb{F}_{N}=N \\
&=&\bigcup_{\gamma \in \Gamma _{A}}(\gamma N
\gamma^{-1})\gamma \mathbb{F}_{A}R \\
&=&\bigcup_{\gamma \in \Gamma _{A}}N\gamma
\mathbb{F}_{A}R\qquad \mathrm{%
because}\,A\,\mathrm{normalizes}\,N \\
&=&N\Gamma _{A}\mathbb{F}_{A}R \\
&=&NAR \\
&=&H\, .
\end{eqnarray*}%
Assume now that
\begin{equation*}
\gamma _{A}\gamma _{N}f_{N}f_{A}r=\sigma _{A}\sigma_{N}g_{N}g_{A}s
\end{equation*}%
for some $\gamma _{A},\sigma _{A}\in \Gamma _{A}$, $\gamma
_{N},\sigma _{N}\in \Gamma _{N}$, $f_{A},g_{A}\in
\mathbb{F}_{A}$,
$f_{N},g_{N}\in \mathbb{F}_{N}$, and $r,s\in R$. Then, as
$A\times
N\times R\ni
(a,n,r)\mapsto anr\in H$ is a diffeomorphism, it follows
that $r=s$. Hence $%
\gamma _{A}\gamma _{N}f_{N}f_{A}=\sigma _{A}\sigma
_{N}g_{N}g_{A}$. But then -- as $A$ normalizes $N$ --
\begin{eqnarray*}
\gamma _{N}f_{N} &=&(\gamma _{A}^{-1}\sigma _{A})\sigma
_{N}g_{N}(g_{A}f_{A}^{-1}) \\
&=&(\gamma _{A}^{-1}\sigma _{A}g_{A}f_{A}^{-1})\left(
(g_{A}f_{A})^{-1}\sigma _{N}g_{N}(g_{A}f_{A}^{-1})\right)
\end{eqnarray*}%
Hence $\gamma _{A}^{-1}\sigma _{A}g_{A}f_{A}^{-1}=1$ or
\begin{equation*}
\gamma _{A}f_{A}=\sigma _{A}g_{A}\,.
\end{equation*}%
As the union $\Gamma _{A}\mathbb{F}_{A}$ is disjoint, it
follows that $%
\gamma _{A}=\sigma _{A}$ and $f_{A}=g_{A}$. But then the
above
implies that
\begin{equation*}
\gamma _{N}f_{N}=\sigma _{N}g_{N}\,.
\end{equation*}%
But then - again because the union is disjoint -- it follows
that
$\gamma _{N}=\sigma _{N}$ and $f_{N}=g_{N}$.
\end{proof}

Our first application of this theorem is to give a simple
proof of
the main result, Theorem 4.2, of \cite{O02}, without using
the
results of \cite{BT96}. We will reformulate those results so
as to
include sampling on irregular grids, see also \cite{A03}.
Let us
first recall some definition before we state the results.

\begin{definition} Let $\cH$ be a separable Hilbert space,
and let $\J$ be a finite our countable infinite index set. A
sequence $\{f_j\}_{j\in\J}$ in $\cH$ is called a
\textit{frame} if
there exists constants $0<A\le B<\infty$ such that for all
$x\in
\cH$ we have.
$$A\|x\|^2\le \sum_{j\in\J}|(x,f_j)|^2\le B\|x\|^2\, .$$
$\{f_j\}_{j\in\J}$ is a \textit{tight frame} if we can
choose
$A=B$ and a \textit{normalized tight frame} or
\textit{Parceval
frame} if we can choose $A=B=1$.
\end{definition}

\begin{example}\label{e-frames}
Let $\cH$ be a Hilbert space and $\cK\subset \cH$ a closed
subspace. Assume that $\{u_n\}$ is a orthonormal basis of
$\cH$.
Let $\mathrm{pr}:\cH\to \cK$ be the orthogonal projection.
Define
$f_j=\pr (u_j)$. Then $\{f_j\}$ is a Parceval frame for
$\cK$. In
fact it is easy to see that every Parceval frame can be
constructed in this way. In particular we can apply this to
the
situation where $(\Omega,\cT)$ is a spectral pair and
$M\subset
\Omega$ is measurable with $|M|>0$. Then
$\{|\Omega|^{-1/2}e_\lambda\}_{\lambda\in \cT}$ is an
orthonormal
basis for $L^2(\Omega)$ and hence $\{f_\lambda :=|\Omega
|^{-1/2}e_\lambda|_{M}\}_{\lambda\in \cT}$ is a Parceval
frame for
$L^2(M)$.
\end{example}

Assume now that $H=ANR$ satisfies the conditions in Lemma
\ref{le-FH}. Let $\Gamma=\Gamma_A\Gamma_N$ and
$\F_H=\F_N\F_A R$
be as in that Lemma. Suppose that $M\subseteq \R^n$ is $H$
invariant and such that there are finitely many open orbits
$U_1,\ldots ,U_k\subseteq M$ such that $|M\setminus
\bigcup_{j=1}^kU_j|=0$. Finally we assume that for each
$\omega_j\in U_j$ the stabilizer of $\omega_j$ in $H$ under
the
action $(h,\omega)\mapsto (h^{-1})^T(\omega)$ is contained
in $R$
and hence compact. Let $\F_j=F_H\cdot \omega_j$ and
$\F=\bigcup_{j=1}^k\F_j$. Then the Lemma \ref{le-FH} implies
that
we have a multiplicative tiling of $M$ as
\begin{equation}\label{eq-Mdisj}
M=\bigcup_{\gamma\in \Gamma }\gamma\cdot \F
=\bigcup_{\gamma\in
\Gamma^{-1}}\gamma^T(\F)
\end{equation}
(up to set of measure zero). If $f:\R^n\to \C$ is a function
and
$a\in \GL (n,\R)$ let $L_af(x):=f(a^{-1}x)$.

\begin{theorem} Let the notation be as above. Suppose that
$\{e_t|_{\F}\}_{t\in\cT}$ is a frame for $L^2(\F)$. Let
$\varphi
=\chi_\F^\vee$. Then the sequence $\{\pi
((t,\gamma)^{-1})\varphi\}_{(t,\gamma)\in\cT\times \Gamma}$
is a
frame for $L^2(M)$ with the same frame bounds. In particular
$\{\pi ((t,\gamma)^{-1})\varphi\}_{(t,\gamma)\in\cT\times
\Gamma}$
is a Parceval frame for $L^2_M(\R^n)$ if and only if
$\{e_t|_{\F}\}_{t\in\cT}$ is a Parceval frame for $L^2(\F)$.
\end{theorem}

\begin{proof} This follows form the fact that by
(\ref{eq-Mdisj}) we have
$$L^2_M(\R^n)\simeq L^2(M)\simeq\bigoplus_{\gamma\in
\Gamma^{-1}}L^2(\gamma^T\F)$$
where the first isomorphism is given by the Fourier
transform.
\end{proof}

Notice, that we can always find as sequence $\{c
e_t|_{\F}\}$,
$c>0$, which is a Parceval frame for $L^2(M)$ by taking a
spectral
pair $(\Omega ,\cT)$, such that $\F\subset \Omega$, i.e., a
parallelepiped $\Omega$.

There are several ways to state different versions of the
above
theorem. In particular one can have different groups
$H_j=A_jN_jR_j$, with compact stabilizers, such that each of
them
has finitely many open orbits, $U_{j,1},\ldots ,U_{j,k_j}$
such
that $\R^n =\bigcup_{j,l} \Gamma_j U_{j,l}$ a disjoint
union. But
we will not state all those obvious generalizations, but
only
notice the following construction from \cite{FO02}. We refer
to
the Appendix for more details. In \cite{FO02} the authors
started
with a prehomogeneous vector  space $(L,V)$ of parabolic
type, see
\cite{BR} for details. Then $L$ has finitely many open
orbits in
$V$, but in general the stabilizer of a point is not
compact. To
deal with that, the authors constructed for each orbits
$U_j$ a
subgroup $H_j=A_jN_jR_j$ such that $U_j$ is up to measure
zero a
disjoint union of open $H_j$ orbits $U_{j,i}$. It turns out,
that
it is not necessary to pick a different group for each
orbit, the
same group $H=H_j$ works for all the orbits.

\begin{theorem}\label{t-foANR} Let $H=ANR$ be one of the
group constructed in
\cite{FO02}. Then $H$ is admissible.
\end{theorem}
\begin{proof} This follows from Theorem 3.6.3 and Corollary
3.6.4
in \cite{BR}.
\end{proof}

\begin{remark} The statement in \cite{BR} is in fact
stronger
than the above remark. In most cases the group $AN$ has
finitely
many open orbits. This group acts freely and is therefore
admissible. The only exception is the so-called Type III
spaces,
where the group $ANR$ has one orbit and is admissible.
\end{remark}

\begin{example}[$\mathbb{R}^{+}\mathrm{SO}(n)$]
Take $A=\R^+\mathrm{id}$, $R=\SO(n)$, the group of
orientation
preserving rotations in $\R^n$, and $N=\{\id\}$. Then
$H=\R^+\SO
(n)$ is the group of dilations and orientation preserving
rotations. Notice that $g^{-1}=g^{T}$ if $g\in
\mathrm{SO}(n)$ and
therefore $g\cdot \omega
=g(\omega )$. The group $H$ has two orbits $\{0\}$ and
$\mathbb{R}%
^{n}\setminus \{0\}$. The stabilizer of $e_{1}=(1,0,\ldots
,0)^{T}$ is isomorphic to $\mathrm{SO}(n-1)$. In particular
the
stabilizer group is compact. It follows that
$\mathbb{R}^{+}\mathrm{SO}(n)$ is admissible. In fact any
function
with compact support in $\mathbb{R}^{n}\setminus \{0\}$ is,
up to
normalization, the Fourier transform of a admissible
wavelet.
\end{example}

\begin{example}[Diagonal matrices]
Let $H$ be the group of diagonal matrices $H=\{d(\lambda
_{1},\ldots ,\lambda _{n})\mid \lambda _{j}\not=0\}$. Thus
$A=\{d(\lambda_1,\ldots ,\lambda_n)\mid \forall j\, :\,
\lambda_j>0\}$ and $R=\{d(\epsilon_1,\ldots ,\epsilon_n)\mid
\epsilon_j=\pm 1\}$. The group $N$ is trivial. Then $H$ has
one
open and dense orbit
\begin{equation*}
U=\{(x_{1},\ldots ,x_{n})^{T}\mid (j=1,\ldots
,n)\,x_{j}\not=0\}=H\cdot (1,\ldots ,1)^{T}\,.
\end{equation*}%
The stabilizer of $(1,\ldots ,1)^{T}$ is trivial and hence
compact. It follows that $H$ is admissible. We can also
replace
$H$ by the connected group $A$. Then we have $2^{n}$ open
orbits
parametrized by $\epsilon \in \{-1,1\}^{2}$
\begin{equation*}
U_{\epsilon }=\{(x_{1},\ldots ,x_{j})^{T}\mid
\mathrm{sign}(x_{j})=\epsilon _{j}\}=H\cdot (\epsilon
_{1},\ldots
,\epsilon _{n})\,.
\end{equation*}%
The stabilizers are still compact and hence $H$ is
admissible.
\end{example}

\begin{example}[Upper triangular matrices]
Let $H$ be the group of upper triangular $2\times
2$-matrices of
determinant $1$,
\begin{equation*}
H=\left\{ \left(
\begin{matrix}
a & t\cr0 & 1/a%
\end{matrix}%
\right) \mid a\not=0,\,t\in \mathbb{R}\right\} \,.
\end{equation*}%
Here $A$ is the group of diagonal matrices with $a>0$, $N$
is the
group up upper triangular matrices with $1$ on the main
diagonal
and $R=\{\pm \id\}$, Then we have one open orbit of full
measure
given by
\begin{equation*}
U=\{(x,y)^{T}\mid y\not=0\}=H\cdot e_{2}
\end{equation*}%
where $e_{2}=(0,1)^{T}$. The stabilizer of $e_{2}$ is
trivial
which implies that $H$ is admissible.
\end{example}

\section{Construction of wavelet sets}

\noindent We apply now our construction in section
\ref{s-wset} to
discrete subgroups of $\GL (n,\R)$. We start by the
following
reformulation of Theorem \ref{t-LWWW} for discrete groups.
Our aim
is later to apply it to the discrete subgroup $\Gamma_A$
from the
last section. As before we use the notation $a\cdot x =
(a^{-1})^T(x)$.

\begin{lemma} \label{LWWWdiscrete} Let $D$ be a discrete
subgroup
$\GL(n, \R)$.  If for almost every $x\in \brn$, there exists
an
$\epsilon > 0$ such that $D^x_\epsilon$ is finite, then
there
exists a measurable function $h$ such that
\begin{equation}
\sum_{d\in D} |h(d^T x)|^2 = 1 \quad a.e.
\end{equation}
\end{lemma}

We have the following improvement of Lemma
\ref{LWWWdiscrete}.

\begin{proposition}\label{LWWWdimprove} Let $D$ be a
discrete subgroup $\GL(n, \R)$.  If for almost every $x\in
\brn$, there exists an $\epsilon > 0$
such that $D^x_\epsilon$ is finite, then there exists a
measurable
function $g$ of the form $g = \chi_K$ such that
\begin{equation}\label{dcc}
\sum_{d\in D} |g (d^Tx)|^2 = 1 \quad a.e.
\end{equation}
\end{proposition}

\begin{proof}  We first recall some notation and preliminary
results from \cite{LWWW2002}.
For an open ball $B\subset \brn$, we define the orbit
density
function $f_B:\brn\to [0,\infty]$ by
\begin{equation}
f_B(x) = \mu(\{d\in D\mid d^Tx \in \overline{B}\}),
\end{equation}
where $\mu$ is counting measure.  Lemma 2.6 of
\cite{LWWW2002}
asserts that
\begin{equation}
\Omega_0 := \{x\in \brn \mid  D^x_\epsilon \,\,{\mathrm
{non-compact}}
\,\,\forall \epsilon > 0\}
 = \{x\in \brn \mid  f_B(x) = 
\infty,\,\, \forall B : B\cap \cO_x \not=\emptyset\}
\end{equation}
Now, let $\cB = \{B_j\}, \,j\in \N,$ be an enumeration of
the
balls in $\brn$ having rational centers and positive
rational
radii.  Let $F_j = F_{B_j}$.  We claim that
\begin{equation}
\brn = \bigcup_{j\ge 1} \{x\in \brn \mid  f_j(x) = 1\} \bigcup
\Omega
\bigcup N,
\end{equation}
where
\begin{equation}
N:=\bigcup_{d\in D, d\not= \id} \{x\in \brn\mid  d^Tx= x\}.
\end{equation}
To see this, suppose that $x\not\in (\Omega_0 \bigcup N)$. 
Then,
there exists an open ball $B$ such that $B\cap \cO_x \not=
\emptyset$, and $f_B(x) < \infty$.  Since $B \cap \cO_x\not
=\emptyset$, there is a $d_0\in D$ such that $d_0^Tx\in B$.
Therefore, there is a $j$ such that $d_0^Tx\in B_j\subset
B$; in
particular $B_j\cap \cO_x \not=\emptyset$ and $\infty >
f_j(x)
> 0$. Now, write $\{d\in D, \, d^Tx \in \overline{B_j}\} =
\{d_0^T,\ldots d_k^T\}$.
Since $x\not\in N$, the $d_i^Tx = d_j^Tx$ only if $i = j$.
Hence,
there is an open set $O$ such that $\mu(\{d\in D \mid  d^Tx \in
O\} =
1$.  Choose $j$ such that $d_0^Tx\in B_j\subset O$ so that
$f_j(x)= 1$.

Continuing along the lines of \cite{LWWW2002}, let
$$\Omega_1 = \{x\in \brn\mid  f_1(x) = 1\}$$
and
$$\Omega_j = \{x\in \brn\mid  f_j(x)
=1\}\setminus(\Omega_1\bigcup\cdots\bigcup\Omega_{j-1})\,
$$ The
sets $\{\Omega_j\}_{j\ge 1}$ form a disjoint collection of
Borel
sets such that $\brn \setminus (\bigcup_{j=1}^\infty
\Omega_j)$
has measure 0 (it is a subset of $\Omega \bigcup N$). Let us
define
\begin{equation}
g(x) = \sum_{j = 1}^\infty \chi_{\Omega_j}(x)
\chi_{\overline
{B_j}}(x)
\end{equation}
and $g(x) = 0$ for $x\not\in (\bigcup_{j\ge 1} \Omega_j)$.
Note
that
\begin{equation}
g(x) = \chi_K\, ,\quad K = \bigcup_{j=1}^\infty (\Omega_j
\cap
\overline{B_j}),
\end{equation}
so all that is needed to complete the proof is to show
$\{d^TK\mid 
d\in D\}$ is a tiling of $\brn$, equivalently, $\sum_{d\in
D}
g(d^Tx) = 1$ a.e.  This is a special case of the argument in
\cite{LWWW2002}, which we outline now.

First, note that if $x\in \brn$ such that $f_j(x) = 1$ for
some
smallest $j$, then there is a unique $d\in D$ such that
$dx\in
\overline {B_j}$.  Since $f_j$ is constant on orbits,
$d^Tx\in
\Omega_j \cap \overline{B_j}$ and $d^Tx \not \in
(\Omega_1\bigcup\cdots \bigcup\Omega_{j-1})$. Therefore, $x
\in d\cdot K$ and $\bigcup_{d\in D} d^T K = \brn$ up to  a set of
measure $0$.

For disjointness, since $d^T \Omega_j  = \Omega_j$, it
suffices to
check that $(\Omega_j \cap \overline{B_j}) \cap (\Omega_j
\cap
\overline{B_j})a$ has measure 0 for all $d$ not the identity
$\id$.  If $x\in (\Omega_j \cap \overline{B_j})$, then
$f_j(x) =
1$.  If, in addition, $x = d^T\omega$ for some $\omega \in
(\Omega_j \cap \overline{B_j})$, then $d^{-1}x\in
\overline{B_j}$
which means that $d\cdot x = x$ and $d = \id$ since $f_j(x)
= 1$.
\end{proof}

\begin{theorem}\label{groupsexist} Let $D$ be a discrete
subgroup $\GL(n, \R)$
that contains an expansive matrix, and $\cL\subset \brn$ a
full-rank lattice.  Then, there exists a $(D, \cL)$ wavelet
set.
\end{theorem}

\begin{proof}  By Proposition \ref{LWWWdimprove}, there
exists a
function $g = \chi_K$ such that equation \ref{dcc} holds.
Therefore, $\{d^TK\mid  d\in D\}$ tiles $\brn$. Thus, by Theorem
\ref{improvement} that there exists a $(D, \cL)$ wavelet
set.
\end{proof}

One final comment is that in all of the above
considerations, the
set $D$ is assumed to invariant under multiplication by an
expansive matrix.  Removing this condition seems to be very
hard.
Indeed, even when the set $D = \{a^j\mid  j\in \Z\}$, it is not
clear
what happens when $a$ is not an expansive matrix.  In this
case,
the interplay between dilations and translations becomes
crucial
in understanding when there exists a wavelet set. For
example, let
$a = \begin{pmatrix}
2 & 0\cr0 & 2/3%
\end{pmatrix}%
$, $D = \{a^j\mid  j\in\mathbb{Z}\}$, and $\mathcal{T} =
\mathbb{Z}^2$. It is easy to see that there is a set of
finite
measure $\Omega$ such that $\{a^j \Omega\}$ tiles
$\mathbb{R}^2$.
However, there exist lattices $\cL_1$ and $\cL_2$ such that
there
are no $(D, \cL_1)$ wavelet sets, yet there are $(D, \cL_2)$
wavelet sets \cite{S}. Hence, in the non-expansive case, it
is not
enough to simply prove the existence of sets that tile via
dilations and translations separately.

We will now apply this to the discrete subgroup
$\Gamma_A\subset
A$ from the last section, where $A$ is as in Theorem
\ref{t-foANR}.

\begin{theorem}\label{t-expinA} Let $H=ANR$ be one of the
group constructed in
\cite{FO02}. Then the following holds:
\item The group $A$ contains a discrete co-compact subgroup
$\Gamma_A\subset
A$ such that $E_A=\{d\in \Gamma_A\mid d\, \,  \mathrm{is\,\,
expansive}\}$ is a non-trivial subsemigroup. In particular
there
exists an expansive matrix $a$ such that $\Gamma_A$ is $a$
invariant.
\item Let $\Gamma_A$ and $\Gamma_A^+$ be as in the Appendix
and let
$\Gamma =\Gamma_A\Gamma_N$. Then $\Gamma \Gamma_A^+\subset
\Gamma$.
\end{theorem}
\begin{proof} (1) follows from Lemma \ref{le-exansive}
in the appendix. For (2) we recall that $\Gamma_A^+$ is
contained
in the center of $ANR$. Hence
\begin{eqnarray*}
\Gamma\Gamma_A^+ &=&\Gamma_A\Gamma_N\Gamma_A^+\cr
&=&\Gamma_A\Gamma_A^+\Gamma_N\cr &\subset
&\Gamma_A\Gamma_N=\Gamma\, .
\end{eqnarray*}
\end{proof}

We have now proved, using Theorem \ref{groupsexist} the
following
theorem:
\begin{theorem} Let the notation be as in Theorem
\ref{t-expinA}. Let
$\cL$ be a full rank lattice in $\R^n$. Then there exists a
$(\Gamma_A,\cL)$ wavelet set and a $(\Gamma,\cL)$ wavelet
set.
\end{theorem}

This Theorem gives several examples of non-groups of
dilations for
which wavelet sets exist. Unfortunately from the point of
view
characterizing sets $(\mathcal{D}, \mathcal{T})$ for which
wavelet
sets exist, if one starts with the set $\mathcal{D}$, one
still
has to rely on the existence of an object external to the
set
$\mathcal{D}$ for the existence of wavelet sets. It would
also be
interesting to remove the condition that $\cL$ is a lattice.

\section{\protect\bigskip Symmetric cones}

\noindent
In this section we discuss the important example of
homogeneous cones in $%
\mathbb{R}^{n}$. Those cones show up in several places in
analysis. As an example one can take Hardy spaces of
holomorphic
function on tube type domains $\R^n+i\oplus \Omega$
\cite{SW}. An
excellent reference for harmonic analysis on symmetric cones
is
the book by J. Faraut and A. Koranyi \cite{FK}. A nonempty
open
subset $\Omega \subset
\mathbb{R}^{n}$ is called an \textit{open (convex) cone} if
$\Omega $ is convex and $%
\mathbb{R}^{+}\Omega \subseteq \Omega $. Let $\Omega $ be an
open
cone,
define the \textit{dual cone} $\Omega ^{\ast }$ by%
\begin{equation*}
\Omega ^{\ast }:=\{v\in \mathbb{R}^{n}\mid \forall u\in
\overline{\Omega }%
\setminus \{0\}~:~(v,u)>0\}~.
\end{equation*}%
If $\Omega ^{\ast }$ is nonempty, then $\Omega ^{\ast }$ is
a open cone. $%
\Omega $ is \textit{self-dual} if $\Omega =\Omega ^{\ast }$.
Let%
\begin{equation*}
\mathrm{GL}(\Omega )=\{g\in \mathrm{GL}(n,\mathbb{R})\mid
g(\Omega
)= \Omega \}~.
\end{equation*}%
Then $\Omega $ is \textit{homogeneous} if
$\mathrm{GL}(\Omega )$
acts transitively on $\Omega $. From now on we assume that
$\Omega
$ is a self-dual homogeneous cone. Let $g\in
\mathrm{GL}(\Omega )$
and $u\in
\overline{\Omega }\setminus \{0\}$. Then $g(u)\in
\overline{\Omega }%
\setminus \{0\}$. Hence if $v\in \Omega =\Omega ^{\ast }$,
then%
\begin{equation*}
(g^{T}(v),u)=(v,g(u))>0~.
\end{equation*}%
Thus $g^{T}(v)\in \Omega ^{\ast }=\Omega $. It follows that
$\mathrm{GL}%
(\Omega )$ is invariant under transposition, and hence
reductive. Let $e\in \Omega $. Then%
\begin{equation*}
K=\mathrm{GL}(\Omega )^{e}=\{g\in \mathrm{GL}(\Omega )\mid
g(e)=e\}~.
\end{equation*}%
Let $\theta (g)=(g^{-1})^{T}$. Then it is always possible to
choice $e$ such that $K=\{g\in \mathrm{GL}(\Omega )\mid
\theta
(g)=g\}=\mathrm{SO}(n)\cap
\mathrm{GL}(n,\mathbb{R})$. Define the \textit{Lie algebra}
of $\mathrm{GL}%
(\Omega )$ by%
\begin{equation*}
\mathfrak{gl}(\Omega ):=\{X\in M(n,\mathbb{R)}\mid \forall
t\in \mathbb{R}%
~:~e^{tX}\in \mathrm{GL}(\Omega )\}~.
\end{equation*}%
Then $\mathfrak{gl}(\Omega )$ is invariant under the Lie
algebra
automorphism $\dot{\theta}(X)=-X^{T}$. Let%
\begin{equation*}
\mathfrak{k}=\{X\in \mathfrak{gl}(\Omega )\mid
\dot{\theta}(X)=X\}
\end{equation*}%
and%
\begin{equation*}
\mathfrak{s}=\{X\in \mathfrak{gl}(\Omega )\mid
\dot{\theta}(X)=-X\}=\mathrm{%
Symm}(n,\mathbb{R})\cap \mathfrak{gl}(\Omega )
\end{equation*}%
where $\mathrm{Symm}(n,\mathbb{R})$ stand for the space of
symmetric matrices. Let $\mathfrak{a}$ be a maximal subspace
in
$\mathfrak{s}$ such
that $[X,Y]=XY-YX=0$ for all $X,Y\in \mathfrak{a}.$ Notice
that $(X,Y)=%
\mathrm{Tr}(XY^{T})$ is an inner product on
$\mathfrak{gl}(\Omega
)$ and
that, with respect to this inner product,
$\mathrm{ad}(X):\mathfrak{gl}%
(\Omega )\rightarrow \mathfrak{gl}(\Omega )$, $Y\mapsto
\lbrack
X,Y]$
satisfies%
\begin{equation*}
\mathrm{ad}(X)^{T}=\mathrm{ad}(X^{T})~.
\end{equation*}%
Hence the algebra $\{\mathrm{ad}(X)\mid X\in \mathfrak{a}\}$
is a
commuting
family of self adjoint operator on the finite dimensional
vector space $%
\mathfrak{gl}(\Omega )$. Hence there exists a basis
$\{X_{j}\}_{j}$ of $%
\mathfrak{gl}(\Omega )$ consisting of joint eigenvectors of
$\{\mathrm{ad}%
(X)\mid X\in \mathfrak{a}\}$. Let
$\mathfrak{z}(\mathfrak{a})$ be
the zero eigenspace, i.e., the maximal subspace of
$\mathfrak{gl}(\Omega )$ commuting with all $X\in
\mathfrak{a}$.
Then there exists a finite subset $\Delta
\subset \mathfrak{a}^{\ast }\setminus \{0\}$ such that with%
\begin{equation*}
\mathfrak{gl}(\Omega )^{\alpha }=\{Y\in \mathfrak{gl}(\Omega
)\mid
\forall X\in \mathfrak{a}~:~\mathrm{ad}(X)Y=\alpha (X)Y\}
\end{equation*}%
we have%
\begin{equation*}
\mathfrak{gl}(\Omega )=\mathfrak{z}(\mathfrak{a})\oplus
\bigoplus_{\alpha \in \Delta }\mathfrak{gl}(\Omega )^{\alpha
}~.
\end{equation*}%
Notice that if $\alpha \in \Delta $ then $-\alpha \in \Delta
$. In
fact, if
\begin{equation}
X\in \mathfrak{gl}(\Omega )^{\alpha }\Longrightarrow
X^{T}\in \mathfrak{gl}%
(\Omega )^{-\alpha }~.  \label{eq-nbar}
\end{equation}%
Let $\mathfrak{a}^{\prime }=\{X\in \mathfrak{a}\mid \forall
\alpha
\in \Delta ~:~\alpha (X)\not=0\}$. Then
$\mathfrak{a}^{\prime }$
is open and dense in $\mathfrak{a}$. In particular
$\mathfrak{a}^{\prime }\not=\{0\}$. Fix $Z\in
\mathfrak{a}^{\prime
}$ and let $\Delta ^{+}=\{\alpha \in \Delta
\mid \alpha (Z)>0\}$. Then $\Delta =\Delta ^{+}\bigcup
-\Delta ^{+}$, and if $%
\alpha ,\beta \in \Delta ^{+}$ are such that $\alpha +\beta
\in
\Delta $,
then $\alpha +\beta \in \Delta ^{+}$. Let%
\begin{equation*}
\mathfrak{n}=\bigoplus_{\alpha \in \Delta
^{+}}\mathfrak{gl}(\Omega )^{\alpha }.
\end{equation*}%
Then $\mathfrak{n}$ is a nilpotent Lie algebra and
$[\mathfrak{a},\mathfrak{n%
}]\subseteq \mathfrak{n}$. In particular it follows that
$\mathfrak{q}=%
\mathfrak{a}\oplus \mathfrak{n}$ is a solvable Lie algebra.
Notice
that the
alebra $\mathfrak{z(a})$ is invariant under transposition.
Hence $\mathfrak{z}%
(\mathfrak{a})=\mathfrak{z}(\mathfrak{a})\cap
\mathfrak{k}\oplus
\mathfrak{a}
$. Because of (\ref{eq-nbar}) it therefore follows that%
\begin{equation*}
\mathfrak{gl}(\Omega )=\mathfrak{k}\oplus \mathfrak{a}\oplus
\mathfrak{n~.}
\end{equation*}%
This decomposition is called \textit{the Iwasawa
decomposition of }$%
\mathfrak{gl}(\Omega )$. Let $A=\{e^{X}\mid X\in
\mathfrak{a}\}$ and $%
N=\{e^{Y}\mid Y\in \mathfrak{n}\}$. Then $A$ and $N$ are Lie
groups, $A$ is abelian, and $aNa^{-1}=N$ for all $a\in A$.
It
follows that $Q:=AN=NA$ is a Lie group with $N$ a normal
subgroup.
Furthermore we have the following \textit{Iwasawa
decomposition of}
$\mathrm{GL}(\Omega )$:

\begin{lemma}[The Iwasawa decomposition]
The map%
\begin{equation*}
A\times N\times K\ni (a,n,k)\mapsto ank\in
\mathrm{GL}(\Omega )
\end{equation*}%
is an analytic diffeomorphism.
\end{lemma}

We note  that the one dimensional group
$Z=\mathbb{R}^{+}\mathrm{id}$ is a subgroup of
$\mathrm{GL}(\Omega
)$ and in fact $Z\subset A$. If $a(\lambda )=\lambda
\mathrm{id}\in Z$, with $\lambda >1$, then $a(\lambda )$ is
expansive. In particular it follows that the set $E$ of
expansive
matrices
in $A$ is a nonempty subsemigroup of $A$. Let
$X_{0}=\mathrm{id}$, $%
X_{1},\ldots ,X_{r}$ be a basis of $\mathfrak{a}$ and let%
\begin{equation*}
\Gamma _{A}=\{\exp (n_{0}X_{0}+\ldots +n_{r}X_{r})\mid
n_{j}\in \mathbb{Z}%
\}~.
\end{equation*}%
Then $A/\Gamma _{A}$ is compact. Furthermore there exists a
discrete subgroup $\Gamma _{N}\subset N$ such that
$N/\Gamma_{N}$
is compact.

Let now $\mathcal{D}=\Gamma $ and $d=\exp (2X_{0})$. Then
$d$ is
expansive and $d\mathcal{D}=\mathcal{D}d\subset
\mathcal{D}$,
because $d$ is central in $\mathrm{GL}(\Omega )$. It follows
that
the results from the previous sections are applicable in
this
case.

\section*{Appendix: Prehomogeneous vector spaces}

\noindent One way to find admissible groups with finitely
many
open orbits is to start with \textit{prehomogeneous vector
spaces}. Those are pairs $(H,V)$ where $H$ is a reductive
Lie
group, say $H^T=H$, and $V$ is a finite dimensional vector
spaces,
such that $H$ has finitely many open orbits in $V$. There is
no
full classification of those spaces at the moment, but a
subclass,
the \textit{prehomogeneous vector spaces of parabolic type}
has
been classified. We refer to \cite{BR} Section 2.11, for
detailed
discussion and references. The problem, from the point of
view of
our work is, that the compact stabilizer condition does not
hold
in general, but as shown in \cite{FO02} one can always
replace $H$
by a subgroup of the form $ANR$ as before, such that $ANR$
is
admissible. Notice that, by using either $ANR$ or
$A^TN^TR^T$,
which satisfies the same conditions, we can consider either
the
standard action on $\R^n$ or the action $(a,x)\mapsto
(a^{-1})^T(x)$. We will use the second action in what
follows.

Let $H=H^T$ be a reductive Lie group acting on
$V=\mathbb{R}^n$.
Then $H$ can be written as $H=LC$ where $C=C^T$ is a vector
group,
isomorphic to a abelian subalgebra $\mathfrak{c}$ of
$\mathfrak{gl} (n,\mathbb{R})=M(n,\mathbb{R})$. The
isomorphism is
simply given by the matrix exponential function
\begin{equation*}
X\mapsto \exp (X)=e^X=\sum_{j=0}^\infty \frac{X^j}{j!}\, .
\end{equation*}
The vectorspace $V$ is graded in the sense that there exists
a
subset $\Delta\subset \mathfrak{c}^*$ such that
\begin{equation}
V=\bigoplus_{\alpha\in\Delta}V_\alpha
\end{equation}
where
\begin{equation}
V_\alpha =\{v\in V\mid (\forall \in \mathfrak{c})\, :\,
X\cdot
v=\alpha (X)v\}\, .
\end{equation}
If $c=\exp (X)\in C$ and $\lambda\in \mathfrak{c}^*$, then
we
write $c^{\lambda}=e^{\lambda(X)}$. In particular $c \cdot
v=c^{\alpha}v$ for all $v\in V_\alpha$.

Denote by $\mathop{\mathrm{pr}} _\alpha$ the projection onto
$V_\alpha$ along $\bigoplus_{\beta\not=\alpha}V_\beta$.

\begin{lemma}
We have $H\cdot V_\alpha\subset V_\alpha$ for all $\alpha
\in
\Delta$. Furthermore if $v\in V$ and $H\cdot v$ is open,
then
$\mathop{\mathrm{pr}} _\alpha(v)\not= 0$ for all $\alpha\in
\Delta$.
\end{lemma}

\begin{proof}
Let $c\in C$, $h\in H$ and $v\in V_\alpha$. As $C$ is
central in
$C$ it follows that $c\cdot (h\cdot v)=h\cdot (c\cdot
v)=c^{\alpha}h\cdot v$.
\end{proof}

The set $\Delta$ has the properties that $0\notin\Delta$, if
$\alpha\in
\Delta$, then $-\alpha\notin\Delta$, and finally there
exists $%
\alpha_1,\ldots ,\alpha_k\in \Delta$ such that if $\alpha\in
\Delta$, then there are $n_1,\ldots ,n_r\in \mathbb{N}_0$
such
that

\begin{equation}  \label{simpleroots}
\alpha =n_1\alpha_1+\ldots +n_r\alpha_r\,.
\end{equation}
For $\alpha\in \Delta$ let $\mathcal{N}_\alpha=\{X\in
\mathfrak{c}\mid \alpha (X)=0\}$. Then $\bigcup_{\alpha\in
\Delta}\mathcal{N}_\alpha$ is a
finite union of hyperplanes and hence the complement is open
and dense in $%
\mathfrak{c}$. Let $\mathfrak{c}^+$ be a connected component
of
the complement of $\bigcup\mathcal{N}_\alpha$. Because of
(\ref{simpleroots}) we can choose $\mathfrak{c}^+$ such that
\begin{equation}  \label{expanding}
\forall X\in \mathfrak{c}^*\forall \alpha \in \Delta\, :\,
\alpha
(X)>0\, .
\end{equation}
Notice that $\mathfrak{c}^+$ is convex, $\mathfrak{c}^+
+\mathfrak{c}%
^+\subset \mathfrak{c}^+$ and
$\mathbb{R}^+\mathfrak{c}^+\subset \mathfrak{c}%
^+$.

\begin{lemma}\label{le-exansive}
The group $A$ contains a non-trivial abelian semigroup of
expanding matrices.
\end{lemma}

\begin{proof}
Let $C^+:=\exp (\mathfrak{c}^+)$. Suppose that $a,b\in C^+$.
Choose $X,Y\in
\mathfrak{c}^+$ such that $a=\exp (X)$ and $b=\exp(Y)$. Then
$%
ab=\exp(X+Y)\in C^+$. Thus $C^+$ is a semigroup. Let
$a=\exp(X)$
be as above. Let $v=\sum_{\alpha}v_\alpha\in V$ with
$v_\alpha\in$, then
\begin{equation*}
\exp (X)\cdot v=\sum_{\alpha}e^{\alpha (X)}v_\alpha
\end{equation*}
and $e^{\alpha (X)}>1$ because $\alpha (X)>0$ for all
$\alpha$.
\end{proof}

Choice a bais $X_1,\ldots ,X_r$ of $\fa$ such that the
vectors
$X_1,\ldots ,X_k$ and the vectors $X_{k+1},\ldots X_r$ form
a
basis for the orthogonal complement of $\fc$ in $\fa$. Here
we use
the inner product $(X,Y)=\mathrm{Tr}(XY^T)$. As $\fc^+$ is
an open
cone in $\fc$ we can chooce $X_j\in \fc^+$, $j=1,\ldots ,k$.
Let
$$\Gamma_A:=\{\exp (\sum_{j=1}^r n_jX_j\mid \forall j\, :\,
n_j\in \Z\}\, . $$
Then $\Gamma_A$ is a cocompact, discrete subgroup of $A$ and
every
element of
$$\Gamma_A^+:=\{\exp (n_1X_1+\ldots + n_kX_k)\mid
\forall j\, :\, n_j\ge 0\}\setminus \{\id\}$$ is expansive.
As
$$\exp (n_1X_1+\ldots + n_rX_r)\exp (m_1X_1+\ldots +
m_rX_r)=
\exp ((n_1+m_1)X_1+\ldots + (n_r+m_r)X_r)$$ it follows that
$\Gamma_A^+$ is a subsemigroup of $\Gamma_A$ such that
$\Gamma_A\Gamma_A^+\subset \Gamma_A$. Thus $\Gamma_A$ is
$\gamma$
invariant for all $\gamma\in \Gamma_A^+$. We have therefoe
proved
the following Lemma:
\begin{lemma}\label{le-expansive} There exists a co-compact
discrete subgroup $\Gamma$ of
$A$ and a subsemigroup $\Gamma_A^+$ such that each element
of
$\Gamma_A^+$ is expansive and $\Gamma_A\Gamma_A^+\subset
\Gamma_A$.
\end{lemma}

\end{document}